\newtheorem{theoreme}{Th\'eor\`eme}
\newtheorem{corollaire}[theoreme]{Corollaire}
\newtheorem{proposition}[theoreme]{Proposition}

\newcommand{\cartier}{{P.Cartier \& D.Foata}}
\newcommand{\viennot}{{X.G.Viennot}}

\newcommand{\blanc}{\hspace{3ex}}
\newcommand{\cqfd}{.\hfill $\Box$ \par}
\newcommand{\exemple}[1]{\bf Exemple #1:\rm}
\newcommand{\exsuite}[1]{\bf Exemple #1 \rm (suite):}
\newcommand{\preuve}{{\em Preuve: }}

\newcommand{\alea}{al\'eatoire}
\newcommand{\carac}{caract\'eristique}
\newcommand{\decomp}{d\'ecomposition}
\newcommand{\dld}{demi-largeur droite}
\newcommand{\eqr}{\'equerre}
\newcommand{\mono}{mono\"\i de}
\newcommand{\Mono}{Mono\"\i de}
\newcommand{\PC}{{\bf pc}$(G)$}
\newcommand{\pcom}{partiellement commutatif}
\newcommand{\pcoms}{partiellement commutatives}
\newcommand{\proba}{probabilit\'e}
\newcommand{\senum}{s\'erie \'enum\'eratrice}
\newcommand{\senums}{s\'eries \'enum\'eratrices}
\newcommand{\sgen}{s\'erie g\'en\'eratrice}
\newcommand{\sgens}{s\'eries g\'en\'eratrices}
\newcommand{\theor}{th\'eor\`eme}

\newcommand{\eq}[1]{(\ref{eq.#1})}
\newcommand{\fig}[1]{figure~\ref{fig.#1}}
\newcommand{\sct}[1]{section~\ref{sec.#1}}
\newcommand{\thref}[1]{\theor~\ref{th.#1}}
\newcommand{\cor}[1]{corollaire~\ref{th.#1}}

\documentclass[a4paper,11pt,french]{article}
\usepackage[french]{babel}
\usepackage{latexsym}
\usepackage{graphicx}

\title{Mod\`eles avec particules dures, animaux dirig\'es,
  et s\'eries en variables \pcoms}
\author{
  J.B\'etr\'ema et J.G.Penaud\\ 
  \\
  UFR Math\'ematiques et Informatique\\
  Universit\'e Bordeaux I\\
  \& LaBRI, URA CNRS 1304
}
\date{Mai 1993}

\begin{document}
\maketitle

\selectlanguage{french}

\begin{abstract}
  Nous exposons syst\'ematiquement les relations entre
  les mod\`eles de gaz avec particules dures,
  l'\'enum\'eration des animaux dirig\'es,
  et l'alg\`ebre des s\'eries formelles en variables \pcoms,
  d'apr\`es les travaux de \viennot.
  Nous donnons une solution com\-pl\`ete simplifi\'ee,
  en utilisant cette alg\`ebre, de l'\'enum\'eration
  des animaux dirig\'es sur r\'eseau plan.
  L'article inclut un programme rapide de g\'en\'eration al\'eatoire
  de ces animaux, avec distribution uniforme, et des
  images obtenues \`a l'aide de ce programme.
\end{abstract}

\def\abstract{
\small 
\begin{center}
{\bf Abstract\vspace{-.5em}\vspace{0pt}} 
\end{center}
\quotation 
}

\begin{abstract}
  We give a systematic presentation of relations between
  lattice gas models with hard-core interactions,
  enumeration of directed-site animals,
  and the algebra of formal power-series in the partially
  commutative case, along the work of \viennot.
  We present a complete and simplified solution,
  using this algebra, of the directed animals enumeration problem
  in two dimensions, including a fast program for random
  generation of such animals (with uniform distribution)
  and images produced by this program.
\end{abstract}

\newpage

\tableofcontents

\newpage

\listoffigures

\section{Introduction}

Certains physiciens (voir par exemple
  \cite{baxter,biggs,dhar1,dhar2,kast})
s'int\'eressent \`a des mod\`eles thermodynamiques
extr\^emement sch\'ematiques, mais o\`u apparaissent d\'ej\`a des
ph\'enom\`enes {\em critiques}, dont ils peuvent faire une \'etude
d\'etaill\'ee. Dans ces mod\`eles, des particules sont plac\'ees
sur certains sommets d'un graphe $G$, appel\'es {\em sites}
(au plus une particule par site);
la disposition des particules forme 
une {\em configuration}, et l'{\em \'energie\/}
d'une configuration est la somme de deux s\'eries de termes:
les premiers repr\'esentent la contribution individuelle
de chaque particule,
les seconds l'{\em interaction\/} entre particules {\em voisines\/}
(c'est \`a dire reli\'ees par une ar\^ete du graphe).

Les configurations \'etudi\'ees sont de tr\`es grande taille,
et on s'int\'eresse \`a leurs propri\'et\'es statistiques;
cette \'etude se ram\`ene \`a celle de la
{\em fonction de partition} du mod\`ele,
not\'ee traditionnellement $Z$,
qui joue un r\^ole voisin d'une {\em fonction g\'en\'eratrice}
en combinatoire; d'ailleurs le calcul de $Z$ se ram\`ene souvent
explicitement \`a celui d'une s\'erie \'enum\'eratrice:
voir en particulier \cite{kast}.

\viennot\ a d\'ecouvert -- voir par exemple \cite{bourbaki} --
que, dans le cas particulier des mod\`eles dits avec
{\em particules dures},
cette \'etude est li\'ee \`a celle du \mono\ \pcom,
introduit en 1969 par \cartier\ \cite{cartier},
et qui fait l'objet aujourd'hui de recherches actives
dans un domaine \`a priori \'eloign\'e,
celui de la mod\'elisation des algorithmes parall\`eles
en informatique th\'eorique.
Cet article expose cette approche
de fa\c con syst\'ematique, et inclut la solution selon
cette ligne, due aux auteurs, du probl\`eme de l'\'enum\'eration
des animaux dirig\'es plans.

La \sct{modeles} pr\'esente la classe
de mod\`eles \'etudi\'es, et contient les d\'efini\-tions issues
de la physique, en particulier celle de la
fonction de partition.

La \sct{series} introduit la
notion  d'{\em empilement} de configurations,
les op\'era\-tions alg\'ebriques sur ces objets,
les th\'eor\`emes fondamentaux de cette alg\`ebre
(qui est celle des s\'eries en variables \pcoms),
et leurs corollaires pour les mod\`eles avec particules dures.

La \sct{animaux} pr\'esente les animaux dirig\'es
sur un r\'eseau, et montre que, sous certaines conditions,
l'\'enum\'eration de ces animaux, sur un r\'eseau de dimension $d$,
se ram\`ene directement \`a celle d'empilements,
sur un r\'eseau de dimension $d-1$.

Dans la \sct{solution}, on ach\`eve les calculs
dans le cas le plus simple, celui du mod\`ele de gaz
sur r\'eseau lin\'eaire, qui correspond aux 
animaux dirig\'es sur r\'eseau plan.
Dans ce cas, les s\'eries en variables \pcoms\ v\'erifient des
\'equations alg\'ebriques simples, qui fournissent \`a la
fois les \senums\ des animaux, et des bijections avec des objets
combinatoires classiques, les chemins (ou mots) de Motzkin.

Ces bijections sont utilis\'ees \sct{aleatoires}, pour programmer
la construction d'animaux dirig\'es plans al\'eatoires,
avec distribution uniforme, en temps lin\'eaire.
La \sct{aleatoires} est illustr\'ee de dessins typiques
d'animaux de grande taille, produits par notre logiciel;
leur aspect filiforme surprend l'intuition.

La \sct{perspectives} \'enum\`ere les principaux probl\`emes
ouverts, en particulier celui d'une interpr\'etation
combinatoire de la solution de R.Baxter pour le
mod\`ele des hexagones durs.

\section{Mod\`eles en m\'ecanique statistique}
\label{sec.modeles}

\subsection{Mod\`ele d'Ising et mod\`ele de gaz}

\label{sec.ising}

Dans le {\em mod\`ele d'Ising\/} du ferromagn\'etisme,
chaque site $i$ est occup\'e par un petit aimant,
avec deux orientations oppos\'ees possibles,
not\'es \mbox{$\sigma_{i} = \pm 1$,}
et l'\'energie est donn\'ee par:
\begin{equation}
  \label{eq.fer}
  E = - H \sum_{{\rm site\ }i} \sigma_{i} -
  J \sum_{{\rm ar\hat ete\ }(i,j)} \sigma_{i}\sigma_{j}
\end{equation}
o\`u $H$ d\'esigne le champ magn\'etique externe,
et $J$ la constante d'interaction.
L'\'energie due au champ externe est minimale lorsque
les aimants sont orient\'es dans le sens du champ;
lorsque la constante $J$ est positive,
l'\'energie d'interaction est minimale lorsque
tous les aimants sont orient\'es dans le m\^eme sens,
d'o\`u le nom de ferromagn\'etisme
---~on ne s'\'etonnera pas de constater que l'\'energie
minimale est n\'egative, car $E$ est d\'efini \`a une
constante additive pr\`es.
Lorsque $J<0$, chaque aimant pousse ses voisins \`a
s'orienter en sens contraire de lui-m\^eme:
on parle d'antiferromagn\'etisme.

Il existe un mod\`ele analogue des gaz:
chaque site $i$ est occup\'e ou vide,
ce qu'on note $s_{i} = 1$ ou $0$,
et l'\'energie vaut:
\begin{equation}
  \label{eq.gaz}
  E = - \mu \sum_{{\rm site\ }i} s_{i} -
  \epsilon \sum_{{\rm ar\hat ete\ }(i,j)} s_{i} s_{j}
\end{equation}
$\mu$ est appel\'e potentiel chimique;
le cas $\epsilon>0$ correspond \`a une attraction entre particules,
et  $\epsilon<0$ \`a une r\'epulsion;
le cas $\epsilon=0$ est celui, sans int\'er\^et ici,
du gaz parfait.

Dans le mod\`ele du ferromagn\'etisme,
il existe une interaction entre chaque paire de sites voisins;
dans le mod\`ele du gaz, cette interaction
n'existe que si les sites sont occup\'es.
Pour comparer les deux mod\`eles, on effectue
le changement de variable $\sigma_{i} = 2 s_{i}-1$; l'\'energie
d'interaction associ\'ee \`a l'ar\^ete $(i,j)$ devient alors:
\[
\sigma_{i}\sigma_{j} = 4 s_{i} s_{j} - 2 s_{i} -2 s_{j} + 1
\]
Lorsqu'on somme cette quantit\'e pour toutes les ar\^etes,
le terme lin\'eaire $- 2 s_{i}$ appara\^\i t avec un coefficient
\'egal au degr\'e $q_{i}$ du site $i$
(le degr\'e d'un sommet est d\'efini comme le nombre d'ar\^etes incidentes).
Donc \eq{fer} devient:
\[
  E = - 2 H \sum_{{\rm site\ }i} s_{i} +
  2 J \sum_{{\rm site\ }i} q_{i} s_{i} -
  4 J \sum_{{\rm ar\hat ete\ }(i,j)} s_{i} s_{j}
\]
en ignorant la constante $nH-mJ$,
o\`u $n$ (resp. $m$) d\'esigne le nombre de sommets
(resp. ar\^etes) du graphe.
Si le graphe est {\em r\'egulier}, chaque degr\'e $q_{i}$ est \'egal \`a
une constante $q$, et on voit que dans ce cas
les deux mod\`eles sont \'equivalents, avec:
\[ \mu = 2 H - 2 q J, \blanc \epsilon = 4 J \]
(cf. \cite[section 1.9]{baxter}).

\subsection{Fonction de partition}

\label{sec.partition}

L'axiome essentiel de la m\'ecanique statistique
(formule de Gibbs),
affirme que la \proba\ d'observation d'une configuration donn\'ee,
d'\'energie totale $E$,
est proportionnelle \`a
\( e^{-\beta E} \),
o\`u $\beta$ est une constante positive.
$\beta$ d\'epend du milieu ext\'erieur, dans lequel est plong\'e le
mod\`ele; plus pr\'eci\-s\'ement:
\[ \beta = \frac{1}{kT} \]
o\`u $k$ d\'esigne la constante de Boltzmann,
et $T$ la {\em temp\'erature} absolue.
Lorsque $\beta$ est grand, les configurations de forte
\'energie deviennent tr\`es improbables:
le milieu ext\'erieur est froid;
inversement, lorsque $\beta$ tend vers z\'ero,
toutes les configurations tendent \`a devenir \'equiprobables:
le milieu ext\'erieur est chaud.

La fonction
\[
  Z(T) = \sum e^{-\beta E} = \sum e^{-E/kT}
\]
o\`u la sommation est \'etendue \`a
toutes les configurations possibles, est appel\'ee
{\em fonction de partition}
---~d\'efinir l'\'energie $E$ \`a une constante additive pr\`es,
a pour effet de d\'efinir $Z$ \`a une constante multiplicative pr\`es.
La \proba\ d'observer une configuration donn\'ee,
d'\'energie $E$, vaut donc:
\begin{equation}
  \label{eq.proba}
  \frac{e^{-\beta E}}{Z(T)}
\end{equation}

La temp\'erature n'est pas en g\'en\'eral la seule variable
dont d\'epende la fonction $Z$: 
d'autres grandeurs physiques, comme le champ magn\'etique,
peuvent intervenir dans le calcul de l'\'energie d'une configuration,
et donc de $Z$. Les physiciens montrent  que l'\'etude statistique
du mod\`ele se ram\`ene toujours \`a l'\'etude de la fonction $Z$
et de ses d\'eriv\'ees ou  d\'eriv\'ees partielles;
la section suivante illustrera ce principe g\'en\'eral. 

\subsection{Mod\`eles avec particules dures}

\label{sec.particules-dures}

Un cas limite des mod\`eles de gaz pr\'ec\'edents
est obtenu en supposant que la r\'epulsion entre particules
est si forte
que les configurations comportant des particules voisines ont une
probabilit\'e n\'egligeable d'exister;
autrement dit, chaque particule est suffisamment dure,
et occupe un volume tel, que
deux sites voisins ne peuvent \^etre occup\'es simultan\'ement.
Les seules configurations possibles correspondent alors aux
ensembles de sommets qu'on appelle {\em stables}
en th\'eorie des graphes:
deux sommets distincts d'une m\^eme configuration ne sont jamais
reli\'es par une ar\^ete.
L'\'energie d'une telle configuration
de $n$ \'el\'ements se r\'eduit \`a $-n \mu$,
et en posant
\[
   t = e^{\beta \mu} = e^{\mu / kT}
\]
la fonction de partition vaut:
\begin{equation}
   \label{eq.Z}
   Z(t) = \sum \alpha_{n} t^{n} 
\end{equation}
o\`u $\alpha_{n}$ d\'esigne le nombre de stables de taille $n$.
En utilisant \eq{proba}, on en d\'eduit que
le nombre moyen de particules vaut:
\begin{equation}
  \label{eq.moyenne}
  \frac{\sum n \alpha_{n} t^{n}}{Z(t)} = 
    \frac{tZ'(t)}{Z(t)} = t \frac{d}{dt}\ln Z(t)
\end{equation}

Cet exemple est caract\'eristique du r\^ole de la fonction de partition
$Z$ et de ses d\'eriv\'ees dans un calcul de moyenne.

\subsection{Limite thermodynamique}

\label{sec.thermo}

La d\'efinition de la fonction de partition pour un graphe infini
n'est pas imm\'ediate:
les formules pr\'ec\'edentes n'ont plus de sens,
puisque le nombre de configurations est infini.
D'autre part la seule classe int\'eressante de graphes infinis,
dans ce contexte, est celle des {\em r\'eseaux}.

En math\'ematiques (voir par exemple \cite{conway}), un r\'eseau
de dimension $d$ d\'esigne habituellement, dans un espace affine,
un ensemble de points isomorphe \`a ${\bf Z}^{d}$;
les points associ\'es, par l'isomorphisme, \`a ceux
dont les coordonn\'ees valent 0 ou 1, forment
le parall\'elotope de base $B$; dans un espace euclidien,
un r\'eseau est cristallographique
si les sym\'etries orthogonales par rapport aux faces
de $B$ laissent le r\'eseau invariant.

\begin{figure}[htbp]
\begin{centering}
\setlength{\unitlength}{0.02cm}
\newcommand{\Sommet}{\circle*{6}}

\begin{picture}(480,240)(80,435)


\put(180,435){\makebox(0,0){(a)}}

\multiput(100,480)(40,0){5}{\Sommet}
\multiput(100,520)(40,0){5}{\Sommet}
\multiput(100,560)(40,0){5}{\Sommet}
\multiput(100,600)(40,0){5}{\Sommet}
\multiput(100,640)(40,0){5}{\Sommet}

\multiput(80,480)(0,40){5}{\line(1,0){200}}
\multiput(100,460)(40,0){5}{\line(0,1){200}}


\put(460,435){\makebox(0,0){(b)}}

\multiput(380,480)(40,0){5}{\Sommet}
\multiput(380,520)(40,0){5}{\Sommet}
\multiput(380,560)(40,0){5}{\Sommet}
\multiput(380,600)(40,0){5}{\Sommet}
\multiput(380,640)(40,0){5}{\Sommet}

\multiput(360,480)(0,40){5}{\line(1,0){200}}
\multiput(380,460)(40,0){5}{\line(0,1){200}}

\put( 360,460){\line( 1, 1){200}}
\put( 360,500){\line( 1, 1){160}}
\put( 360,540){\line( 1, 1){120}}
\put( 360,580){\line( 1, 1){80}}
\put( 400,460){\line( 1, 1){160}}
\put( 440,460){\line( 1, 1){120}}
\put( 480,460){\line( 1, 1){80}}

\put( 560,460){\line( -1, 1){200}}
\put( 560,500){\line( -1, 1){160}}
\put( 560,540){\line( -1, 1){120}}
\put( 560,580){\line( -1, 1){80}}
\put( 520,460){\line( -1, 1){160}}
\put( 480,460){\line( -1, 1){120}}
\put( 440,460){\line( -1, 1){80}}

\end{picture}
\caption{R\'eseaux carr\'es (a) {\em nn} (b) {\em nnn}}
\label{fig.nn}
\end{centering}
\end{figure}
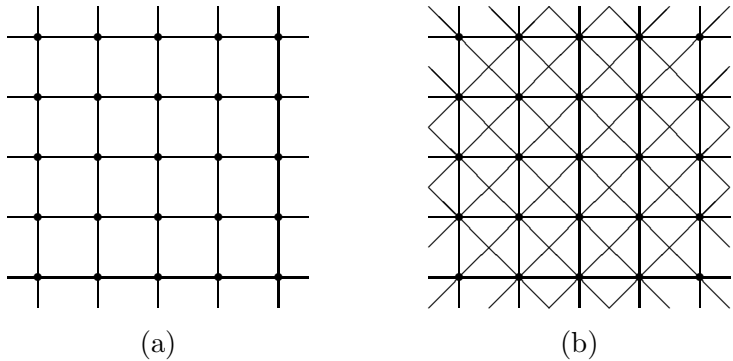

Dans le pr\'esent contexte, un r\'eseau d\'esigne un graphe,
dont l'ensemble des sommets est un r\'eseau cristallographique;
le choix le plus simple pour les ar\^etes est de ne relier
que les voisins les plus proches (r\'eseau {\em nn},
pour ``nearest neighbour'', chez les anglo-saxons),
mais on peut ajouter des ar\^etes qui relient les voisins 
``au second rang''
(r\'eseau {\em nnn}, pour ``next nearest neighbour''), etc\ldots
A un m\^eme r\'eseau cristallographique peuvent donc correspondre
plusieurs graphes:
sur la \fig{nn}-(a),
chaque sommet est de degr\'e 4,
tandis que sur la \fig{nn}-(b),
chaque sommet est de degr\'e 8.

Notons qu'on peut \'etendre la notion de r\'eseau
en affaiblissant les hypoth\`eses sur l'ensemble des sites:
par exemple, si on pave le plan par des hexagones,
les sommets ne forment plus un ensemble isomorphe \`a
${\bf Z}^{2}$ (il manque les centres des hexagones),
mais le graphe obtenu est r\'egulier,
et d\'efinit un mod\`ele d'interaction
int\'eressant, o\`u chaque site est de degr\'e 3.

Pour d\'efinir la fonction de partition d'un r\'eseau $G$,
on consid\`ere l'intersection $G_{k}$ de $G$ avec un
disque de rayon donn\'e, o\`u $k$ d\'esigne
le nombre de sites de $G_{k}$.
Soit $Z_{k}$ la fonction de partition de  $G_{k}$;
$Z_{k}$ n'est pas une quantit\'e proportionnelle
\`a la taille $k$ du syst\`eme
(le nombre de configurations cro\^\i t exponentiellement
avec la taille),
mais la m\'ecanique statistique enseigne que
c'est le cas de $\ln Z_{k}$,
lorsque $k$ est grand;
plus pr\'ecis\'ement, la fonction $Z$ d\'efinie par:
\begin{equation}
   \label{eq.limZ}
   \ln Z = \lim_{k \rightarrow \infty} \frac{\ln Z_{k}}{k}
\end{equation}
existe, et cette limite est ind\'ependante
de la d\'efinition pr\'ecise de $G_{k}$.
Les physiciens appellent ce type de limite la {\em limite thermodynamique};
consulter R.Baxter \cite{baxter}, ou un cours de m\'ecanique statistique,
comme \cite{chandler},
pour une pr\'esentation compl\`ete.

Les ph\'enom\`enes critiques correspondent 
aux {\em singularit\'es} de la fonction de partition $Z$
d\'efinie par \eq{limZ}, et \`a celles de ses
d\'eriv\'ees ou  d\'eriv\'ees partielles.
Le nombre moyen de particules {\em par site},
appel\'e {\em densit\'e moyenne par site},
est donn\'e, sur un r\'eseau infini,
par le membre droit
de la formule \eq{moyenne}
o\`u $Z$ d\'esigne cette fois la fonction de partition pour la
limite thermodynamique.

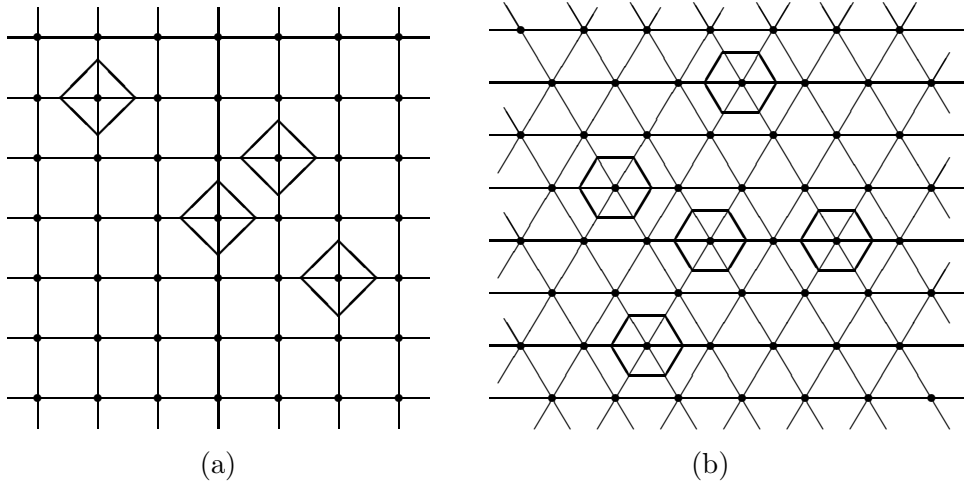
\begin{figure}[htbp]
\begin{centering}
\setlength{\unitlength}{0.02cm}
\newcommand{\Sommet}{\circle*{6}}

\begin{picture}(635,320)(80,435)


\thicklines
\put( 195,600){\line( 1, 1){25}}
\put( 195,600){\line( 1, -1){25}}
\put( 245,600){\line( -1, 1){25}}
\put( 245,600){\line( -1, -1){25}}

\put( 235,640){\line( 1, 1){25}}
\put( 235,640){\line( 1, -1){25}}
\put( 285,640){\line( -1, 1){25}}
\put( 285,640){\line( -1, -1){25}}

\put( 115,680){\line( 1, 1){25}}
\put( 115,680){\line( 1, -1){25}}
\put( 165,680){\line( -1, 1){25}}
\put( 165,680){\line( -1, -1){25}}

\put( 275,560){\line( 1, 1){25}}
\put( 275,560){\line( 1, -1){25}}
\put( 325,560){\line( -1, 1){25}}
\put( 325,560){\line( -1, -1){25}}
\thinlines

\multiput(100,480)(40,0){7}{\Sommet}
\multiput(100,520)(40,0){7}{\Sommet}
\multiput(100,560)(40,0){7}{\Sommet}
\multiput(100,600)(40,0){7}{\Sommet}
\multiput(100,640)(40,0){7}{\Sommet}
\multiput(100,680)(40,0){7}{\Sommet}
\multiput(100,720)(40,0){7}{\Sommet}

\multiput(80,480)(0,40){7}{\line(1,0){280}}
\multiput(100,460)(40,0){7}{\line(0,1){280}}

\put(220,435){\makebox(0,0){(a)}}


\thicklines
\put( 535,565){\line( 1, 0){24}}
\put( 535,605){\line( 1, 0){24}}
\put( 523,585){\line( 3, 5){12}}
\put( 523,585){\line( 3, -5){12}}
\put( 571,585){\line( -3, 5){12}}
\put( 571,585){\line( -3, -5){12}}

\put( 493,495){\line( 1, 0){24}}
\put( 493,535){\line( 1, 0){24}}
\put( 481,515){\line( 3, 5){12}}
\put( 481,515){\line( 3, -5){12}}
\put( 529,515){\line( -3, 5){12}}
\put( 529,515){\line( -3, -5){12}}

\put( 472,600){\line( 1, 0){24}}
\put( 472,640){\line( 1, 0){24}}
\put( 460,620){\line( 3, 5){12}}
\put( 460,620){\line( 3, -5){12}}
\put( 508,620){\line( -3, 5){12}}
\put( 508,620){\line( -3, -5){12}}

\put( 619,565){\line( 1, 0){24}}
\put( 619,605){\line( 1, 0){24}}
\put( 607,585){\line( 3, 5){12}}
\put( 607,585){\line( 3, -5){12}}
\put( 655,585){\line( -3, 5){12}}
\put( 655,585){\line( -3, -5){12}}

\put( 555,670){\line( 1, 0){24}}
\put( 555,710){\line( 1, 0){24}}
\put( 543,690){\line( 3, 5){12}}
\put( 543,690){\line( 3, -5){12}}
\put( 591,690){\line( -3, 5){12}}
\put( 591,690){\line( -3, -5){12}}
\thinlines

\put(430,460){\line( 3, 5){ 170}}
\put(472,460){\line( 3, 5){ 170}}
\put(514,460){\line( 3, 5){ 170}}
\put(556,460){\line( 3, 5){ 150}}

\put(598,460){\line( 3, 5){ 108}}
\put(640,460){\line( 3, 5){ 66}}

\put(406,490){\line( 3, 5){ 150}}
\put(406,560){\line( 3, 5){ 108}}
\put(406,630){\line( 3, 5){ 66}}

\put(580,460){\line( -3, 5){ 170}}
\put(622,460){\line( -3, 5){ 170}}
\put(664,460){\line( -3, 5){ 170}}
\put(706,460){\line( -3, 5){ 170}}

\put(706,530){\line( -3, 5){ 128}}
\put(706,600){\line( -3, 5){ 86}}
\put(706,670){\line( -3, 5){ 44}}

\put(538,460){\line( -3, 5){ 128}}
\put(496,460){\line( -3, 5){ 86}}
\put(454,460){\line( -3, 5){ 44}}

\multiput(400,480)(0,35){8}{\line(1,0){ 315}}

\multiput(442,480)(42,0){7}{\Sommet}
\multiput(421,515)(42,0){7}{\Sommet}
\multiput(442,550)(42,0){7}{\Sommet}
\multiput(421,585)(42,0){7}{\Sommet}
\multiput(442,620)(42,0){7}{\Sommet}
\multiput(421,655)(42,0){7}{\Sommet}
\multiput(442,690)(42,0){7}{\Sommet}
\multiput(421,725)(42,0){7}{\Sommet}

\put(547,435){\makebox(0,0){(b)}}

\end{picture}
\caption{Exemples de configurations}
  dans le mod\`ele des
  (a) carr\'es durs (b) hexagones durs \\
\label{fig.carhex}
\end{centering}
\end{figure}

Le mod\`ele des {\em carr\'es\/} (resp. {\em hexagones\/})
{\em durs\/} correspond \`a l'\'etude d'un gaz,
dont les particules sont dispos\'ees sur un r\'eseau carr\'e
(resp. hexagonal) infini, de telle sorte
que deux sites voisins ne peuvent \^etre simultan\'ement occup\'es.
Pour rendre visible cette derni\`ere condition,
on peut repr\'esenter les particules par des carr\'es
(resp. hexagones) qui ne se chevauchent pas, comme sur la
\fig{carhex}. Dans la suite de cet article, un site occup\'e
par une particule sera simplement indiqu\'e par un disque noir.

Le calcul de la fonction de partition
du mod\`ele des hexagones durs a \'et\'e r\'eussi par Baxter en 1980:
voir \cite{baxter} pour une pr\'esentation compl\`ete; c'est un calcul
complexe, qui fait intervenir entre autres
les identit\'es de Rogers-Ramanujan;
le calcul de la fonction de partition
du mod\`ele des carr\'es durs reste un probl\`eme ouvert.

\newcommand{\astar}{$A^{\ast}$}
\newcommand{\Gbar}{\overline{\Gamma}}
\newcommand{\Tbar}{\overline{\Theta}}
\newcommand{\Pbar}{\overline{\Pi}}

\section{S\'eries en variables \pcoms}

\label{sec.series}

\subsection{Empilements}

\label{sec.empilements}

Soit un {\em graphe} $G$; les sommets de $G$ formeront dans la section suivante
les lettres d'un {\em alphabet}, et leur ensemble sera donc not\'e $A$;
une {\em ar\^ete} de $G$ est une paire (non orient\'ee)
$\{a, b\}$ de sommets distincts.

Deux sommets de $G$ sont dits {\em voisins} s'ils sont identiques, ou reli\'es
par une ar\^ete de $G$; le voisinage $V(a)$ d'un sommet $a$ est l'ensemble
des sommets voisins de $a$; le voisinage $V(B)$ d'une partie
$B \subseteq A$ est d\'efini comme d'habitude par:

\[ V(B) = \bigcup_{b \in B} V(b) \]

Dans toute la suite nous conviendrons qu'une {\em configuration}
$C$ est un ensemble de sommets tel
qu'aucune paire de sommets distincts de $C$ ne soient voisins
(c'est ce qu'on appelle un {\em stable} en th\'eorie des graphes).
Un {\em empilement}
est une suite finie $(C_1, C_2, \ldots, C_n)$ de
configurations non vides telle que:
\begin{equation}
  \label{eq.empilement}
  C_i \subseteq V (C_{i-1}), \blanc 1 < i \leq n 
\end{equation}
L'entier $n$ est la {\em hauteur} de l'empilement;
chaque configuration $C_i$ est appel\'ee 
{\em couche} num\'ero $i$,
et la premi\`ere couche, $C_1$, est la {\em base} de l'empilement.
Si $n=0$, l'empilement est vide, et ne poss\`ede pas de base.

Un empilement peut aussi \^etre d\'ecrit comme un ensemble de
{\em cellules} $(a,i)$, o\`u $a$ est un sommet de $G$, et $i$ un entier
inf\'erieur ou \'egal \`a la hauteur $n$ de l'empilement; $i$ est
la hauteur de la cellule, et la couche num\'ero $i$ est form\'ee
des cellules de hauteur $i$;
les cellules ayant m\^eme projection $a$ sur $A$
forment la {\em fibre} $F_a$ au-dessus de $a$.
On dira pour abr\'eger que deux cellules $(a,i)$ et $(b,j)$ sont voisines si
les sommets $a$ et $b$ sont voisins;
on peut alors reformuler la d\'efinition
d'un empilement de la mani\`ere suivante:
\begin{itemize}
  \item
  deux  cellules distinctes de la m\^eme couche
  ne sont jamais voisines
  \refstepcounter{equation}
  \hfill (\arabic{equation})
  
  \label{eq.cell1}
  
  \item
  pour $i > 1$,
  chaque cellule de la couche $i$
  poss\`ede au moins une voisine
  dans la couche $i - 1$
  \refstepcounter{equation}
  \hfill (\arabic{equation})

  \label{eq.cell2}
\end{itemize}

La {\em taille} $|E|$ d'un empilement $E$ est \'egale, par d\'efinition,
au nombre de cellules qui composent $E$. Une configuration
peut \^etre confondue avec un empilement de hauteur 1;
\`a la configuration vide, correspond l'empilement vide, de taille et de
hauteur nulles.

\begin{quote}
\exemple{1} la \fig{empilement} montre un empilement de hauteur 4,
et de taille 8:
\begin{figure}[htbp]
  \begin{center}
  \setlength{\unitlength}{0.02cm}
  \newcommand{\Noir}{\circle*{10}}

\begin{picture}(575,340)(65,420)
\thinlines
\put(420,480){\line( 1, 0){ 80}}
\put(500,480){\line( 2, 1){ 40}}
\put(540,500){\line(-1, 0){ 80}}
\put(460,500){\line(-2,-1){ 40}}
\put(500,460){\line( 2, 1){120}}
\put(620,520){\line(-1, 0){160}}
\put(460,520){\line(-2,-1){120}}
\put(340,460){\line( 1, 0){160}}
\put(460,520){\line( 0,-1){ 20}}
\put(340,460){\line( 4, 1){ 80}}
\put(500,480){\line( 0,-1){ 20}}
\put(540,500){\line( 4, 1){ 80}}
\put(420,560){\line( 1, 0){ 80}}
\put(500,560){\line( 2, 1){ 40}}
\put(540,580){\line(-1, 0){ 80}}
\put(460,580){\line(-2,-1){ 40}}
\put(500,540){\line( 2, 1){120}}
\put(620,600){\line(-1, 0){160}}
\put(460,600){\line(-2,-1){120}}
\put(340,540){\line( 1, 0){160}}
\put(460,600){\line( 0,-1){ 20}}
\put(340,540){\line( 4, 1){ 80}}
\put(500,560){\line( 0,-1){ 20}}
\put(540,580){\line( 4, 1){ 80}}
\put(420,640){\line( 1, 0){ 80}}
\put(500,640){\line( 2, 1){ 40}}
\put(540,660){\line(-1, 0){ 80}}
\put(460,660){\line(-2,-1){ 40}}
\put(500,620){\line( 2, 1){120}}
\put(620,680){\line(-1, 0){160}}
\put(460,680){\line(-2,-1){120}}
\put(340,620){\line( 1, 0){160}}
\put(460,680){\line( 0,-1){ 20}}
\put(340,620){\line( 4, 1){ 80}}
\put(500,640){\line( 0,-1){ 20}}
\put(540,660){\line( 4, 1){ 80}}
\put(420,720){\line( 1, 0){ 80}}
\put(500,720){\line( 2, 1){ 40}}
\put(540,740){\line(-1, 0){ 80}}
\put(460,740){\line(-2,-1){ 40}}
\put(500,700){\line( 2, 1){120}}
\put(620,760){\line(-1, 0){160}}
\put(460,760){\line(-2,-1){120}}
\put(340,700){\line( 1, 0){160}}
\put(460,760){\line( 0,-1){ 20}}
\put(340,700){\line( 4, 1){ 80}}
\put(500,720){\line( 0,-1){ 20}}
\put(540,740){\line( 4, 1){ 80}}
\put( 80,660){\line( 0,-1){160}}
\put( 80,500){\line( 1, 0){160}}
\put(240,500){\line( 0, 1){160}}
\put(240,660){\line(-1, 0){160}}
\put(120,620){\line( 0,-1){ 80}}
\put(120,540){\line( 1, 0){ 80}}
\put(200,540){\line( 0, 1){ 80}}
\put(200,620){\line(-1, 0){ 80}}
\put( 80,660){\line( 1,-1){ 40}}
\put(200,620){\line( 1, 1){ 40}}
\put(200,540){\line( 1,-1){ 40}}
\put( 80,500){\line( 1, 1){ 40}}
\put(130,600){\makebox(0,0)[lb]{h}}
\put(130,550){\makebox(0,0)[lb]{e}}
\put(185,600){\makebox(0,0)[lb]{g}}
\put(185,550){\makebox(0,0)[lb]{f}}
\put(340,460){\Noir}
\put(620,520){\Noir}
\put(500,540){\Noir}
\put(420,560){\Noir}
\put(540,580){\Noir}
\put(420,640){\Noir}
\put(340,700){\Noir}
\put(500,720){\Noir}
\put(480,420){\makebox(0,0)[lb]{E}}
\put(640,480){\makebox(0,0)[lb]{$C_1$}}
\put(640,560){\makebox(0,0)[lb]{$C_2$}}
\put(640,640){\makebox(0,0)[lb]{$C_3$}}
\put(640,720){\makebox(0,0)[lb]{$C_4$}}
\put(160,460){\makebox(0,0)[lb]{G}}
\put( 60,485){\makebox(0,0)[lb]{a}}
\put(250,485){\makebox(0,0)[lb]{b}}
\put(250,665){\makebox(0,0)[lb]{c}}
\put( 60,665){\makebox(0,0)[lb]{d}}
\end{picture}
  \end{center}
  \caption{Un empilement $E$ de configurations sur le graphe $G$}
  \label{fig.empilement}
\end{figure}
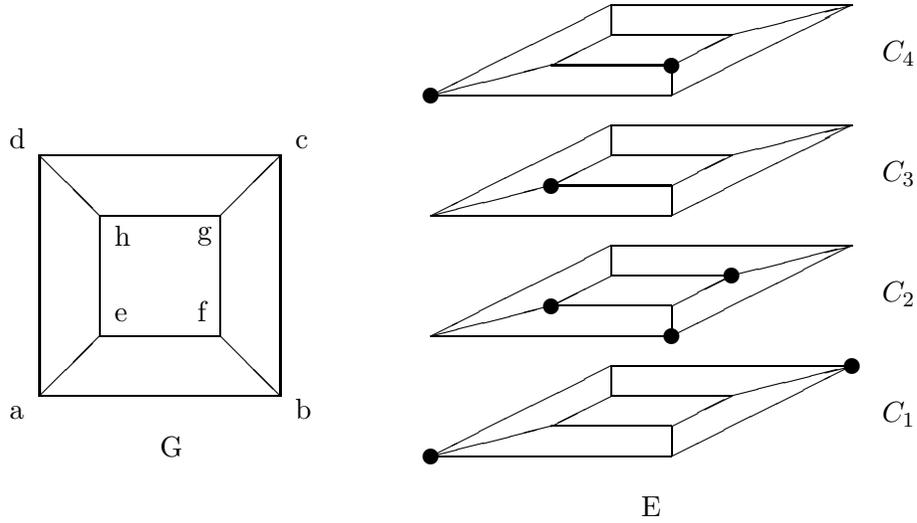

\[
\begin{array}{ll}
C_1 = \{a, c\}  & V(C_1) = \{a,b,c,d,e,g\} \\
C_2 = \{b,e,g\} & V(C_2) = \{a,b,c,e,f,g,h\} \\
C_3 = \{e\}     & V(C_3) = \{a,e,f,h\} \\
C_4 = \{a, f\}
\end{array}
\]
$C_1,C_2,C_3,C_4$ sont des configurations qui v\'erifient la condition
\eq{empilement}, donc
\[
E=( C_1,C_2,C_3,C_4 )
\]
est un empilement, constitu\'e des cellules:
\[ (a,1),(c,1),(b,2),(e,2),(g,2),(e,3),(a,4),(f,4). \]
\end{quote}

La notion d'empilement est due \`a \viennot\ \cite{heaps},
qui repr\'esente les cellules par des {\em pi\`eces} g\'eom\'etriques,
dessin\'ees de telle sorte que deux pi\`eces aient une intersection
non vide si et seulement si elles sont voisines.
Les pi\`eces de la couche $C_1$ sont pos\'ees sur le sol,
celles de la couche $C_2$ sont pos\'ees par dessus,
et ainsi de suite (deux pi\`eces de la m\^eme fibre ont m\^eme
projection sur le sol). Gr\^ace \`a la condition \eq{cell2},
chaque pi\`ece d'un empilement repose sur sa (ou ses) voisine(s)
de la couche inf\'erieure.

\subsection{\Mono\ \pcom}

\label{sec.monoide-pc}

Soit \astar\ le \mono\ libre sur l'alphabet $A$ des sommets du graphe $G$;
le \mono\ \pcom\ associ\'e \`a $G$, not\'e \PC, est le quotient
de \astar\ par les relations de commutation:
\[ ab=ba \]
pour tous sommets $a,b$ non voisins dans $G$ (on dit que $G$ est le
graphe de non-commu\-ta\-tion). Les \'el\'ements de \PC\ sont 
appel\'es {\em traces}, et la trace du mot $u$ est not\'ee $[u]$.
Deux mots $u,v$ de \astar\ ont m\^eme trace si l'on peut passer de l'un \`a
l'autre par une suite de {\em commutations autoris\'ees},
c'est \`a dire d'\'echanges entre lettres cons\'ecutives d\'esignant
des sommets de $G$ non voisins; dans ce cas on dit qu'ils sont
\'equivalents, et on note: $u \sim v$.

A une configuration $C = \{a,b,c,\ldots\}$,
associons la trace du mot $u(C)=abc\ldots$;
par d\'efinition d'une configuration,
les lettres $a,b,c,\ldots$ commutent toutes entre elles,
et la trace $[u(C)]$ est donc bien d\'efinie, car elle
ne d\'epend pas de l'ordre dans lequel on
\'enum\`ere les \'el\'ements de $C$.

A un empilement $E = ( C_1, C_2, \ldots, C_n )$,
associons le produit des $[u(C_i)]$ dans \PC\ 
(si l'empilement est vide, on lui associe la trace du mot vide),
et notons $\Phi$ cette application.

\begin{quote}
\exsuite{1}
la trace $\Phi (E)$ de l'empilement $E$ de la \fig{empilement},
est la classe du mot $acbegeaf$ modulo
les 16 commutations:
$ac=ca, af=fa, \ldots, fh=hf$.
\end{quote}

R\'eciproquement, \`a tout mot $u=a_1 \ldots a_n$ de \astar,
associons l'empilement $E(u)$ constitu\'e des cellules
$( a_i, h(a_i) )$, o\`u la hauteur $h(a_i)$
de la lettre $a_i$
est d\'efinie par r\'ecurrence comme suit:
\begin{equation}
  \label{eq.hauteur}
  h(a_i) = 1 + \max \{h(a_j)\ |\ j < i {\rm \ et\ }
  a_j \in V(a_i) \}
\end{equation}
avec la convention $\max (\emptyset) = 0$.

On v\'erifie facilement que $E(u)$ est bien un empilement --
c'est \`a dire satisfait aux conditions \eq{cell1} et \eq{cell2} --
et que cet empilement ne change pas si l'on applique \`a $u$
une commutation autoris\'ee.
Par passage au quotient, on obtient donc une application $\Psi$ de \PC\
dans l'ensemble des empilements,
et on a le r\'esultat suivant, \'enonc\'e dans \cite{heaps}:

\begin{theoreme}
\label{th.emp-pc}
Les applications $\Phi$ et $\Psi$ sont des bijections,
inverses l'une de l'autre,
entre l'ensemble des empilements sur un graphe $G$ et le \mono\ \pcom\ \PC.
\end{theoreme}
\preuve lorsqu'on associe un mot $u$ \`a un empilement $E$,
en \'enum\'erant les couches de $E$, puis qu'on
calcule les hauteurs des lettres de $u$ par la formule \eq{hauteur},
on retrouve les num\'eros des couches. 
Inversement, lorsqu'on part d'un mot $u$,
la d\'efinition \eq{hauteur} implique que les
lettres de hauteur 1 commutent avec toutes celles qui les pr\'ec\`edent,
et donc $u \sim u_1v$, o\`u $u_1$ contient toutes les
lettres de hauteur 1; et ainsi de suite, ce qui prouve que $u$
est \'equivalent au mot obtenu en \'enum\'erant
les couches de $E(u)$\cqfd

\bigskip

En particulier, les empilements forment un {\em \mono\ }
isomorphe \`a \PC;
notons que le produit des empilements
\(
( C_1, \ldots, C_n )
{\rm \ et\ }
( C'_1, \ldots, C'_p )
\)
n'est pas en g\'en\'eral \'egal \`a la suite de configurations
\( ( C_1, \ldots, C_n, C'_1, \ldots, C'_p ) \),
qui ne constitue pas un empilement, sauf dans le cas exceptionnel o\`u
\( C'_1 \subseteq V(C_n) \).
Dans le langage
g\'eom\'etrique de \viennot, on compose deux empilements
en pla\c cant le second au dessus du premier, et en laissant les
pi\`eces du second tomber sur celles du premier.

\begin{quote}
\exsuite{1}
si on ajoute le sommet $c$ \`a
l'empilement $E$ de la \fig{empilement},
c'est \`a dire si on calcule le produit $E'=E.c$,
on a:
\[
  E'= (C_1,C_2,C'_3,C_4), \blanc C'_3=C_3 \cup \{c\}
\]
car $c$ commute avec tous les \'el\'ements de $C_3$ et
$C_4$, mais pas avec $b \in C_2$.
\end{quote}

Notons aussi que lorsqu'on lit un mot $u$ de droite \`a gauche,
on obtient un nouvel empilement, {\em dual} de celui obtenu
en lisant normalement le mot de gauche \`a droite;
cette transformation correspond g\'eom\'etriquement
\`a une inversion du sens de la gravit\'e.

Gr\^ace au \thref{emp-pc},
les couches d'un empilement peuvent \^etre interpr\'et\'ees
comme une {\em factorisation canonique} d'un \'el\'ement
du \mono\ \pcom; c'est sous cette forme que ce th\'eor\`eme
est apparu pour la premi\`ere fois,
dans \cite[\theor\ 1.2]{cartier}, sous le nom
de {\em V-factorisation}.

\subsection{S\'eries formelles et th\'eor\`eme d'inversion}

\label{sec.inversion}

A partir du \mono\ \PC\
-- que nous ne distinguerons plus d\'esormais du \mono\
des empilements --
on peut construire l'alg\`ebre des s\'eries formelles \`a coefficients
rationnels (par exemple) sur ce \mono:
c'est l'alg\`ebre des s\'eries en variables \pcoms.
Un \'el\'ement de cette alg\`ebre est une application
$\Delta$ qui associe \`a chaque empilement $E$
(ou \`a sa trace)
un coefficient $\Delta(E)$ dans le corps de base choisi;
on note:
\[
  \Delta = \sum_{{\rm empilement\ }E} \Delta(E) E
\]
et le produit est d\'efini naturellement par:
\[
  \left( \sum_{E_1} \Delta_1(E_1) E_1 \right)
  \left( \sum_{E_2} \Delta_2(E_2) E_2 \right)
  =
  \sum_{E} \left(
    \sum_{E_1E_2=E} \Delta_1(E_1) \Delta_2(E_2)
  \right ) E
\]
On v\'erifie facilement qu'il n'existe qu'un nombre fini de
factorisations de $E$ sous la forme $E_1E_2$,
ce qui assure la coh\'erence de la d\'efinition du produit.

Consid\'erons en particulier les s\'eries  suivantes,
qui sont en fait des polyn\^omes lorsque le graphe $G$ est fini:
\begin{equation}
  \label{eq.gamma}
  \Gamma = \sum_{{\rm configuration\ }C}C
  \blanc\blanc
  \Gbar = \sum_{{\rm configuration\ }C}(-1)^{|C|}C
\end{equation}
et les s\'eries analogues (s\'erie \carac\ et s\'erie altern\'ee)
pour les empilements:
\begin{equation}
  \label{eq.theta}
  \Theta = \sum_{{\rm empilement\ }E}E
  \blanc\blanc
  \Tbar = \sum_{{\rm empilement\ }E}(-1)^{|E|}E
\end{equation}
en rappelant que $|E|$ d\'esigne le nombre de cellules de $E$.

\begin{quote}
\exemple{2}
\blanc soit $G$ le graphe:
\begin{picture}(140,40)(20,5)
  \thinlines
  \put(50,10){\line(1,0){100}}
  \put(50,10){\circle*{5}}
  \put(100,10){\circle*{5}}
  \put(150,10){\circle*{5}}
  \put(50,0){\makebox(0,0){$a$}}
  \put(100,0){\makebox(0,0){$b$}}
  \put(150,0){\makebox(0,0){$c$}}
\end{picture}

\medskip\noindent Omettons les crochets [\ ]
pour d\'esigner les traces des mots dans \PC;
le polyn\^ome $\Gbar$ vaut:
\[ 1-a-b-c+ac \]
tandis que la s\'erie $\Theta$ d\'ebute par:
\[ 1+a+b+c+a^2+ab+ac+ba+b^2+bc+cb+c^2+\ldots \]
(noter l'absence du terme $ca$, car $ac=ca$; un comptage
manuel montre que la s\'erie comporte ensuite 21 termes de degr\'e~3).

\end{quote}

On a le r\'esultat fondamental suivant:

\begin{theoreme}
\label{th.inversion}
La s\'erie \carac\ des empilements,  $\Theta$, et
la s\'erie altern\'ee des configurations, $\Gbar$,
d\'efinies par \eq{gamma} et \eq{theta},
sont inverses l'une de l'autre. De m\^eme,
$\Tbar$ est l'inverse de $\Gamma$.
\end{theoreme}

\preuve
soit une configuration C et un empilement E;
l'\'equation $CX=E$, dans le \mono\ des empilements,
a une solution (unique) $X$ si et seulement si:
\( C \subseteq C_1 \),
o\`u $C_1$ d\'esigne la base de $E$.
Donc en effectuant le produit  $\Gbar\:\Theta$, le coefficient d'un empilement donn\'e $E$ non vide vaut:
$\sum (-1)^{|C|}$,
la sommation \'etant \'etendue \`a tous les sous-ensembles $C$ de
$C_1$. Or cette somme est nulle
(c'est la base des arguments combinatoires du type ``inclusion-exclusion''),
ce qui prouve que $\Gbar\:\Theta = 1$.

Le m\^eme raisonnement est valable pour l'\'equation: $XC=E$,
en lisant les mots de droite \`a gauche dans \PC\
(cf. empilement dual)
ce qui prouve que $\Theta\:\Gbar = 1$.
En changeant chaque variable en son oppos\'ee, on obtient aussi:
$\Gamma\:\Tbar = \Tbar\:\Gamma = 1$\cqfd

\bigskip

Ce \theor\ appara\^\i t dans \cite{heaps}, formul\'e
avec des valuations;
sa source, cit\'ee dans \cite{heaps},
est \`a nouveau \cite[\theor\ 2.4]{cartier},
o\`u il est pr\'esent\'e comme un calcul de {\em fonction de M\"obius}
(cf. la preuve).
Cas particuliers:
\begin{itemize}
\item
si $G$ est le graphe {\em complet} sur $A$, on retrouve la formule:
\[
\frac{1}{1-A} = A^{\ast}
\]
dans l'alg\`ebre des s\'eries en variables non commutatives.

\item
si $G$ est le graphe {\em totalement d\'econnect\'e}
--~c'est \`a dire sans ar\^ete~-- sur $A$,
toutes les variables commutent,
la s\'erie $\Gbar$ se factorise en:
\( \prod_{a \in A}(1-a) \)
et la s\'erie $\Theta$ en:
\( \prod_{a \in A}(1+a+a^2+\ldots) \)

\end{itemize}

\subsection{Pyramides et th\'eor\`eme de d\'ecomposition}

\label{sec.pyramides}

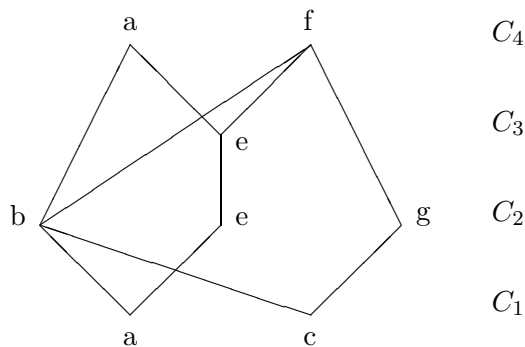
\begin{figure}[htbp]
  \begin{center}
  \setlength{\unitlength}{0.02cm}
  \begin{picture}(315,219)(65,300)
\thinlines
\put(140,320){\line(-1, 1){ 60}}
\put(140,320){\line( 1, 1){ 60}}
\put(260,320){\line(-3, 1){180}}
\put(260,320){\line( 1, 1){ 60}}
\put(200,380){\line( 0, 1){ 60}}
\put(200,440){\line(-1, 1){ 60}}
\put(200,440){\line( 1, 1){ 60}}
\put( 80,380){\line( 3, 2){180}}
\put( 80,380){\line( 1, 2){ 60}}
\put(260,500){\line( 1,-2){ 60}}

\put(380,320){\makebox(0,0)[lb]{$C_1$}}
\put(380,380){\makebox(0,0)[lb]{$C_2$}}
\put(380,440){\makebox(0,0)[lb]{$C_3$}}
\put(380,500){\makebox(0,0)[lb]{$C_4$}}

\put(135,300){\makebox(0,0)[lb]{a}}
\put(255,300){\makebox(0,0)[lb]{c}}

\put( 60,380){\makebox(0,0)[lb]{b}}
\put(210,380){\makebox(0,0)[lb]{e}}
\put(330,380){\makebox(0,0)[lb]{g}}

\put(210,430){\makebox(0,0)[lb]{e}}

\put(135,510){\makebox(0,0)[lb]{a}}
\put(255,510){\makebox(0,0)[lb]{f}}

\end{picture}
  \end{center}
  \caption{Diagramme de Hasse
    de l'empilement $E$ de la \protect\fig{empilement}.}
  \label{fig.hasse}
\end{figure}

A un empilement $E$ est associ\'e un {\em ordre} sur les cellules
de la mani\`ere suivante: on dit que la cellule
$(a,i)$ pr\'ec\`ede $(b,j)$ dans $E$ si:
\begin{quote}
 $a$ et $b$ sont voisins, et $i<j$
\end{quote}
et on consid\`ere l'ordre engendr\'e par cette relation.
Le {\em diagramme de Hasse} (voir \fig{hasse})
rend cet ordre visible:
c'est le graphe dont les sommets sont les \'el\'ements
de l'ensemble partiellement ordonn\'e consid\'er\'e,
et les arcs les couples $(x,y)$ tels que:
\begin{quote}
  $x<y$  et il n'existe pas d'\'el\'ement $z$ tel que $x<z<y$
\end{quote}

Une {\em pyramide} est un empilement dont la base $C_1$ est un singleton.
Consid\'erons un empilement quelconque $E$,
et une cellule $\gamma$ de $E$;
l'ensemble des cellules sup\'erieures ou \'egales \`a $\gamma$,
pour l'ordre ci-dessus, forme une pyramide de base $\gamma$,
qu'on dit {\em engendr\'ee}  par $\gamma$,
et la d\'ecomposition:
\begin{equation}
  \label{eq.decomposition}
  E=XP, \blanc P {\rm \ pyramide\ de\ base\ }\gamma
\end{equation}
est unique.
Si $E$ est repr\'esent\'e par un mot $u$,
la cellule $\gamma$ correspond \`a une occurrence d'une lettre $a$
dans le mot $u$, et la pyramide $P$ au mot form\'e par les
lettres situ\'ees \`a droite de $a$ dans $u$,
et qu'aucune suite de commutations autoris\'ees ne peut faire
glisser \`a gauche de $a$.

\begin{quote}
\exsuite{1}
l'empilement $E$ est repr\'esent\'e par le mot $u\ =\ acbegeaf$,
et poss\`ede les 8 d\'ecompositions suivantes
(une par cellule), de type \eq{decomposition}, d\'eduites de la
\fig{hasse}:

\(
\begin{array}{ll}
  u\,& =\,cg.abeeaf\,=\,aee.cbgaf\,=\,acege.baf\,=\,acbg.eeaf\\
     & =\,acbeea.gf\,=\,acbeg.eaf\,=\,acbegef.a\,=\,acbegea.f
\end{array}
\)
\end{quote}

Aux s\'eries \carac s et altern\'ees des configurations et empilements,
d\'efinies par \eq{gamma} et \eq{theta},
ajoutons celles des pyramides:
\[
  \Pi = \sum_{{\rm pyramide\ }P}P
  \blanc\blanc
  \Pbar = \sum_{{\rm pyramide\ }P}(-1)^{|P|}P
\]
et les s\'eries {\em d\'eriv\'ees}:
\begin{equation}
  \label{eq.derivees}
  \Gamma\,' = \sum_{{\rm configuration\ }C} |C|\ C
  \blanc\blanc
  \Theta' = \sum_{{\rm empilement\ }E} |E|\ E
\end{equation}
(rappelons que $|E|$ d\'esigne le nombre de cellules de $E$).

Notons qu'il n'y a pas de pyramide vide, donc la s\'erie $\Pi$ ne comporte
pas de terme constant.

\begin{quote}
 \exsuite{2} Soit $G$ le graphe de la \sct{inversion};
  la s\'erie $\Pi$ d\'ebute par: 
  \[ a+b+c+a^2+ab+ba+b^2+bc+cb+c^2+\ldots \]
  (suivis de 18 termes de degr\'e~3).
\end{quote}
\begin{theoreme}
  \label{th.pyramides}
  La s\'erie \carac\ des pyramides est la
  ``d\'eriv\'ee logarithmique droite''
  de celle des empilements; plus pr\'ecis\'ement, dans l'anneau
  des s\'eries formelles en variables \pcoms, on a:
  \begin{equation}
    \label{eq.derivee-log-theta}
    \Theta' = \Theta\:\Pi
  \end{equation}
  De m\^eme, la s\'erie altern\'ee des pyramides est
  l'oppos\'ee de la
  ``d\'eriv\'ee logarithmique gauche''
  de la s\'erie \carac\  des configurations:
  \begin{equation}
    \label{eq.derivee-log-gamma}
    \Gamma\,' = - \Pbar\:\Gamma
  \end{equation}
\end{theoreme}
\preuve
l'\'equation \eq{derivee-log-theta}
r\'esulte directement de la d\'ecomposition \eq{decomposition},
car, pour un empilement donn\'e $E$,
il existe $|E|$ d\'ecompositions de ce type
(une par cellule).
L'\'equation \eq{derivee-log-gamma} est une cons\'equence alg\'ebrique de
\eq{derivee-log-theta} et du \theor\ d'inversion (\thref{inversion}), car
la d\'eriv\'ee logarithmique de l'inverse de $\Theta$ est, bien entendu,
l'oppos\'ee de la d\'eriv\'ee logarithmique de $\Theta$.

Pr\'ecis\'ement, en multipliant \`a gauche les deux termes de
\eq{derivee-log-theta}
par l'inverse de $\Theta$, on obtient:
\[
\Gbar\:\Theta' = \Pi
\]
D'autre part, de l'\'egalit\'e:
\[
|E_1 E_2|E_1 E_2 = (|E_1|+| E_2|)E_1 E_2 = (|E_1|E_1) E_2 + E_1(| E_2| E_2)
\]
on d\'eduit par lin\'earit\'e que la formule habituelle de d\'erivation
d'un produit est applicable aux s\'eries en variables \pcoms,
et donc:
\[
(\Gbar\:\Theta)' = \Gbar\,'\:\Theta + \Gbar\:\Theta'
\]
et la formule $\Gbar\:\Theta = 1$ implique donc:
\[
\Gbar\,'\:\Theta = - \Gbar\:\Theta'
\blanc {\rm d'o\grave u} \blanc
\Gbar\,'\:\Theta = - \Pi
\]
On en d\'eduit \eq{derivee-log-gamma} en multipliant \`a droite par $\Gbar$
et en changeant les variables en leurs oppos\'ees,
ce qui \'echange s\'eries \carac s et altern\'ees\cqfd

\subsection{Applications aux mod\`eles de gaz}
\label{sec.applications-gaz}

\begin{corollaire}
  \label{th.moyenne}
  Soit un mod\`ele de gaz avec particules dures, sur un graphe fini $G$,
  et soit:
  \[
     t = e^{\mu / kT}
  \]
  o\`u $\mu$ d\'esigne le potentiel chimique du mod\`ele,
  $k$ la constante de Boltzmann,
  et $T$ la temp\'erature.
   Le nombre moyen de particules vaut:
  \[
     \sum_{n=1}^{\infty}(-1)^{n-1}p_n\,t^n
  \]
  o\`u $p_n$ d\'esigne le nombre de pyramides de taille $n$ sur $G$.
\end{corollaire}

\preuve
Pour toute s\'erie en variables \pcoms\ $\Delta$,
appelons {\em projection} de $\Delta$ la s\'erie ordinaire
$\delta(t)$ obtenue en rempla\c cant
les lettres de l'alphabet $A$ par une m\^eme variable t.
La projection est donc l'application:
\[
  \Delta = \sum_{E} \Delta(E) E \blanc \rightarrow \blanc
  \delta(t) = \sum_{n=0}^{\infty}
    \left( \sum_{|E|=n} \Delta(E) \right) t^n
\]
Comme par hypoth\`ese $G$ est fini, il n'existe qu'un nombre fini
d'empilements de taille $n$, et la s\'erie $\delta$
est bien d\'efinie.

Par d\'efinition, d\'eriver $\Delta$ -- cf. \eq{derivees} --
consiste \`a multiplier chaque coefficient $\Delta(E)$
par la taille $|E|$ de l'empilement; donc la projection de $\Delta'$
se d\'eduit de $\delta$ en rempla\c cant $t^n$ par $nt^n$,
autrement dit vaut $t\,\delta'(t)$.

Les projections de
$\Gamma, \Theta$ et $\Pi$ sont
les s\'eries g\'en\'eratrices
$\gamma(t), \theta(t)$ et $\pi(t)$
qui \'enum\`erent respectivement
les configurations, empilements et pyramides selon leur taille;
la formule \eq{derivee-log-gamma} implique:
\[
  t\,\gamma'(t) = - \pi(-t)\,\gamma (t)
\]

Or la fonction de partition $Z(t)$ du mod\`ele
est \'egale \`a $\gamma(t)$, selon la formule \eq{Z},
\sct{particules-dures}; et le nombre moyen de
particules vaut $t\,Z'(t) / Z(t)$, selon \eq{moyenne}:
le corollaire suit\cqfd

\begin{corollaire}
  \label{th.densite}
  Soit un mod\`ele de gaz avec particules dures, sur un r\'eseau $G$,
  et soit $O$ un sommet de $G$.
  La densit\'e moyenne par site vaut:
  \[
     \sum_{n=1}^{\infty}(-1)^{n-1}p_n\,t^n
  \]
  o\`u $p_n$ d\'esigne le nombre de pyramides de taille $n$ sur $G$,
  de base $O$.
\end{corollaire}

\preuve
sur un r\'eseau, infini par d\'efinition,
projeter une s\'erie n'a en g\'en\'eral plus de sens,
car il existe une infinit\'e d'empilements de taille donn\'ee;
par contre on peut projeter la s\'erie $\Pi_{O}$
des pyramides de base $O$,
et la s\'erie g\'en\'eratrice correspondante $\pi(t)$
ne d\'epend pas du choix de $O$,
gr\^ace \`a l'homog\'en\'eit\'e d'un r\'eseau.
D'autre part, soit  $G_{k}$ une portion finie du r\'eseau,
comportant $k$ sommets
(cf. \sct{thermo}); on se doute que si $k$ est
suffisamment grand, et  $G_{k}$ d\'efinie de fa\c con raisonnable,
le nombre de pyramides de taille $n$ sur  $G_{k}$, et de base $a$,
d\'epend peu du choix de $a$; autrement dit la s\'erie
 $\pi_{k}(t)$ qui \'enum\`ere {\em toutes} les pyramides sur  $G_{k}$,
vaut approximativement $k\,\pi(t)$.
On peut montrer qu'on a exactement:
\begin{equation}
  \label{eq.limite-pi}
   \pi(t) = \lim_{k \rightarrow \infty} \frac{\pi_{k}(t)}{k}
\end{equation}
Par d\'efinition, la densit\'e par site est \'egale
au nombre de particules divis\'e par le nombre de sommets;
le \cor{densite} est donc une cons\'equence du
pr\'ec\'edent\cqfd

\bigskip

Ce corollaire peut \^etre interpr\'et\'e comme la limite
thermodynamique du \cor{moyenne}, au sens de la \sct{thermo};
de m\^eme que les physiciens doivent d\'efinir la fonction de partition
d'un gaz sur un r\'eseau par une moyenne logarithmique,
ici il faut remplacer l'\'enum\'eration des configurations
par celle des pyramides de base fix\'ee, qui seule garde un sens.
Notons cependant que notre d\'emarche,
en dehors de la formule \eq{limite-pi},
est purement alg\'ebrique, et repose sur les deux th\'eor\`emes
fondamentaux de l'alg\`ebre des s\'eries formelles en variables
\pcoms\ (inversion et d\'ecomposition).

\subsection{Empilements stricts}

\label{sec.stricts}

Un empilement est dit {\em strict} si deux couches cons\'ecutives
ne contiennent jamais la m\^eme lettre, autrement dit si la condition
\eq{empilement} est remplac\'ee par:
\[
  C_i \subseteq V (C_{i-1}) - C_{i-1}, \blanc 1 < i \leq n 
\]
Si $u$ est un mot de \astar\ qui repr\'esente l'empilement $E$,
$E$ est strict si et seulement si, pour toute d\'ecomposition:
\[ u = u_1\ a\ u_2\ a\ u_3 \]
il existe une lettre de $u_2$ qui ne commute pas avec $a$
(voir par exemple la formule \eq{hauteur} pour le prouver);
on dira qu'un tel mot est lui aussi strict.

\begin{theoreme}
  \label{th.stricts}
  La s\'erie \carac\ des empilements,  $\Theta$,
  se d\'eduit de celle des empilements stricts en substituant
  \[ \frac{a}{1-a} = a+a^2+a^3+\ldots \]
  \`a chaque lettre $a$ de l'alphabet A. Inversement, on passe
  des empilements aux empilements stricts par la substitution:
  \[ a \rightarrow \frac{a}{1+a} = a-a^2+a^3-\ldots \]
\end{theoreme}
\preuve
la seconde affirmation est une cons\'equence imm\'ediate
de la premi\`ere,
car l'inverse de la fonction $a \rightarrow a/(1-a)$
est la fonction $a \rightarrow a/(1+a)$.

Consid\'erons ensuite un empilement $E$,
et $u$ un mot qui repr\'esente $E$;
\'ecrivons $u$ en regroupant les lettres cons\'ecutives \'egales:
\[
  u = a_1^{r_1} a_2^{r_2} \ldots  a_m^{r_m}
  \ {\rm avec} \
  a_i \neq a_{i+1}
  \ {\rm pour} \
  1 \leq i < m
\]
et parmi tous les mots $u$ \'equivalents qui repr\'esentent
le m\^eme empilement $E$, choisissons-en un avec le nombre $m$
de facteurs minimal. Alors $v = a_1 a_2 \ldots  a_m$
est strict; sinon il existe $i<j$
avec $a_i=a_j$ qui commute avec tout $a_{k}$,
pour $i<k<j$, et on peut diminuer le nombre de facteurs de $u$
en regroupant $ a_i^{r_i}$ et $a_j^{r_j}$
en $a_i^{r_i + r_j}$.

Il reste \`a montrer que la trace du mot strict ainsi obtenu est unique;
soient $v = a_1 a_2 \ldots  a_m$ et
$v' = b_1 b_2 \ldots  b_n$ deux mots stricts
tels que 
\(  u = a_1^{r_1} a_2^{r_2} \ldots  a_m^{r_m} \)
et
\(  u' = b_1^{s_1} b_2^{s_2} \ldots  b_n^{s_n} \)
soient \'equivalents; soit $b_i$
la premi\`ere occurrence de $a_1$ dans $v'$; puisque $u \sim u'$,
$b_i$ commute avec $b_{k}$ pour $k<i$.
Posons $w = b_1 \ldots b_{i-1} b_{i+1} \ldots b_n$;
on a donc $v' \sim a_1 w$.
D'autre part, si par exemple
$r_1 > s_i$, il existe $j>i$ tel que $b_j = b_i$
commute avec $b_{k}$ pour $i<k<j$,
ce qui contredit le fait que $v'$ soit strict;
on fait un raisonnement sym\'etrique si $r_1 < s_i$.
Donc $r_1 = s_i$, et on peut simplifier l'\'equivalence
$u \sim u'$ par $a_1^{r_1}$;
par r\'ecurrence sur la longueur des mots, on sait que
$a_2 \ldots  a_m \sim w$, d'o\`u $v \sim v'$\cqfd

\section{Animaux dirig\'es}

\label{sec.animaux}

\subsection{D\'efinition}

Les physiciens appellent {\em animal} un ensemble de points,
ou {\em cellules}, sur un r\'eseau, qui se d\'eveloppe \`a
partir d'une configuration initiale par adjonction de cellules
sur les sites vides voisins de ceux occup\'es
(deux sites sont voisins s'ils sont reli\'es par une ar\^ete
du r\'eseau).
La {\em taille} d'un animal est \'egale au nombre de ses cellules.
 Si la configuration initiale est connexe,
en particulier si elle est r\'eduite
\`a un seul point, l'animal est connexe:
deux cellules distinctes
peuvent toujours \^etre reli\'ees,
dans le r\'eseau, par un chemin dont
tous les sommets appartiennent \`a l'animal.
La g\'en\'eration de ces animaux ob\'eit \`a diff\'erentes lois,
selon qu'ils servent \`a mod\'eliser des probl\`emes de percolation,
de croissance de cristaux ou de polym\`eres, etc\ldots

\begin{figure}[htbp]
\setlength{\unitlength}{0.015cm}
\newcommand{\Noir}{\circle*{14}}

\begin{picture}(800,290)(100,390)
\thinlines


\put(220,420){\Noir}
\multiput(220,460)(40,0){5}{\Noir}
\put(220,500){\Noir}
\multiput(380,500)(40,0){3}{\Noir}
\multiput(220,540)(40,0){7}{\Noir}
\put(220,580){\Noir}
\put(300,580){\Noir}
\multiput(460,580)(40,0){2}{\Noir}
\multiput(220,620)(40,0){3}{\Noir}
\multiput(300,660)(40,0){6}{\Noir}

\multiput(220,420)(0,40){7}{\line(1,0){280}}
\multiput(220,420)(40,0){8}{\line(0,1){240}}

\put(90,530){\makebox(0,0)[lb]{Directions}}
\put(90,500){\makebox(0,0)[lb]{privil\'egi\'ees}}
\put(140,460){\vector( 0, 1){ 40}}
\put(140,460){\vector( 1, 0){ 40}}

\put(90,415){\makebox(0,0)[lb]{Source}}
\put(170,420){\vector( 1, 0){ 30}}

\put(340,380){\makebox(0,0)[lb]{(a)}}


\multiput(680,420)(-24,40){7}{\line(1,0){240}}
\multiput(680,420)(48,0){6}{\line(-3, 5){144}}
\put(680,420){\line( 3, 5){120}}
\put(656,460){\line( 3, 5){120}}
\put(632,500){\line( 3, 5){ 96}}
\put(608,540){\line( 3, 5){ 72}}
\put(584,580){\line( 3, 5){ 48}}
\put(560,620){\line( 3, 5){ 24}}
\put(728,420){\line( 3, 5){ 96}}
\put(776,420){\line( 3, 5){ 72}}
\put(824,420){\line( 3, 5){ 48}}
\put(872,420){\line( 3, 5){ 24}}

\put(680,420){\Noir}
\multiput(704,460)(48,0){4}{\Noir}
\put(680,500){\Noir}
\put(728,500){\Noir}
\put(824,500){\Noir}
\multiput(656,540)(96,0){3}{\Noir}
\multiput(680,580)(48,0){2}{\Noir}
\put(704,620){\Noir}
\multiput(680,660)(48,0){2}{\Noir}

\put(770,380){\makebox(0,0)[lb]{(b)}}

\put(560,460){\vector( 1, 0){ 48}}
\put(560,460){\vector(-2, 3){ 24}}
\put(560,460){\vector( 2, 3){ 24}}

\put(550,415){\makebox(0,0)[lb]{Source}}
\put(625,420){\vector( 1, 0){ 30}}

\end{picture}
\caption{Animaux dirig\'es sur r\'eseaux (a) carr\'e (b) hexagonal}
\label{fig.animaux}
\end{figure}
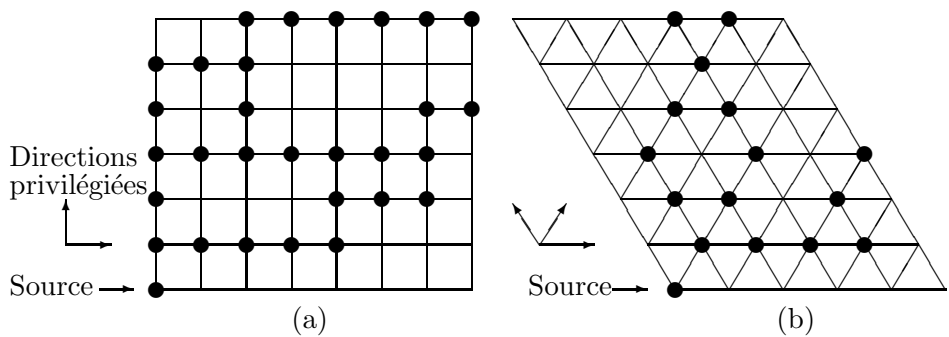

Les {\em animaux dirig\'es} sont contraints de se d\'evelopper
dans des directions privil\'egi\'ees, d\'efinies en orientant
les directions des ar\^etes du r\'eseau;
pour chaque cellule il existe donc dans l'animal un chemin
orient\'e depuis la configuration initiale, appel\'ee {\em source},
vers la cellule;
on suppose en g\'en\'eral que la source est r\'eduite \`a un point.
La \fig{animaux} 
montre des exemples d'animaux dirig\'es sur r\'eseaux plans;
notons que le cas (b) peut \^etre aussi obtenu \`a partir d'un
r\'eseau carr\'e {\em nnn} (voir \fig{nn},
\sct{thermo} ), avec les directions
privil\'egi\'ees Est, Nord et Nord-Est;
un tel r\'eseau est souvent appel\'e {\em triangulaire}.

Nous allons montrer, \sct{equivalence}, que sur certains
types de r\'eseaux, les animaux dirig\'es correspondent \`a
des empilements. Cette \'equivalence passe par le {\em coloriage}
des r\'eseaux, qui permet de construire une variante
des empilements, et de d\'efinir clairement les r\'eseaux
pour lesquels l'\'equivalence est observable; c'est l'objet
des deux sections suivantes. 

\subsection{Empilements colori\'es}
\label{sec.empilements-colories}

Nous avons pr\'esent\'e, \sct{monoide-pc},
la bijection fondamentale entre mots en variables
partiellement commutatives,
et empilements; la clef de cette bijection
r\'eside dans la r\'epartition en couches
des lettres successives d'un mot, de telle sorte que
l'ordre des couches sp\'ecifie enti\`erement l'ordre
des lettres dans le mot, aux commutations autoris\'ees pr\`es.
Autrement dit, on ajoute un nouveau sommet $a$ \`a un empilement
en respectant la r\`egle fondamentale suivante:
la hauteur de $a$ est strictement sup\'erieure \`a celles
de toutes les fibres $F_b$, pour $b$ voisin de $a$.

Il est possible d'introduire des variantes dans la d\'efinition des
couches, tout en conservant cette r\`egle;
on obtient alors de nouvelles vari\'et\'es d'empilements,
toutes \'equivalentes entre elles, puisque isomorphes
au \mono\ \pcom. La variante {\em standard},
\'etudi\'ee jusqu'\`a pr\'esent, consiste \`a
placer $a$ aussi bas que possible -- voir formule \eq{hauteur}.
Dans cette section, nous ajouterons une contrainte,
qui impose \`a tous les sommets d'une couche d'\^etre
de la m\^eme couleur.

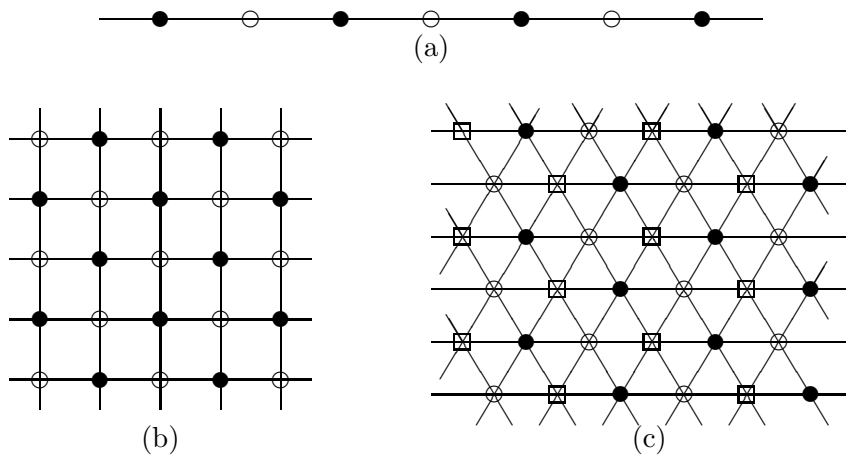
\begin{figure}[htbp]
\begin{centering}
\setlength{\unitlength}{0.02cm}
\newcommand{\Noir}{\circle*{10}}
\newcommand{\Blanc}{\circle{10}}
\newcommand{\Carre}{
  \makebox(0,0)[r]{\framebox(10,10){}}
}

\begin{picture}(635,310)(80,510)
\thinlines


\put(180,800){\line( 1, 0){440}}
\multiput(220,800)(120,0){4}{\Noir}
\multiput(280,800)(120,0){3}{\Blanc}
\put(400,780){\makebox(0,0){(a)}}


\put(180,560){\Noir}
\put(260,560){\Noir}
\put(140,560){\Blanc}
\put(220,560){\Blanc}
\put(300,560){\Blanc}

\put(180,600){\Blanc}
\put(260,600){\Blanc}
\put(140,600){\Noir}
\put(220,600){\Noir}
\put(300,600){\Noir}

\put(180,640){\Noir}
\put(260,640){\Noir}
\put(140,640){\Blanc}
\put(220,640){\Blanc}
\put(300,640){\Blanc}

\put(180,680){\Blanc}
\put(260,680){\Blanc}
\put(140,680){\Noir}
\put(220,680){\Noir}
\put(300,680){\Noir}

\put(180,720){\Noir}
\put(260,720){\Noir}
\put(140,720){\Blanc}
\put(220,720){\Blanc}
\put(300,720){\Blanc}

\multiput(120,560)(0,40){5}{\line(1,0){200}}
\multiput(140,540)(40,0){5}{\line(0,1){200}}
\put(220,520){\makebox(0,0){(b)}}


\put(430,530){\line( 3, 5){ 128}}
\put(472,530){\line( 3, 5){ 128}}
\put(514,530){\line( 3, 5){ 128}}
\put(556,530){\line( 3, 5){ 107}}
\put(598,530){\line( 3, 5){ 65}}

\put(406,560){\line( 3, 5){ 108}}
\put(406,630){\line( 3, 5){ 66}}

\put(538,530){\line( -3, 5){ 128}}
\put(580,530){\line( -3, 5){ 128}}
\put(622,530){\line( -3, 5){ 128}}
\put(664,530){\line( -3, 5){ 128}}
\put(664,600){\line( -3, 5){ 86}}
\put(664,670){\line( -3, 5){ 44}}

\put(496,530){\line( -3, 5){ 86}}
\put(454,530){\line( -3, 5){ 44}}

\put(400,550){\line( 1, 0){ 275}}
\put(442,550){\Blanc}
\put(482,550){\Carre}
\put(526,550){\Noir}
\put(568,550){\Blanc}
\put(608,550){\Carre}
\put(652,550){\Noir}

\put(400,585){\line( 1, 0){ 275}}
\put(419,585){\Carre}
\put(463,585){\Noir}
\put(505,585){\Blanc}
\put(545,585){\Carre}
\put(589,585){\Noir}
\put(631,585){\Blanc}

\put(400,620){\line( 1, 0){ 275}}
\put(442,620){\Blanc}
\put(482,620){\Carre}
\put(526,620){\Noir}
\put(568,620){\Blanc}
\put(608,620){\Carre}
\put(652,620){\Noir}

\put(400,655){\line( 1, 0){ 275}}
\put(419,655){\Carre}
\put(463,655){\Noir}
\put(505,655){\Blanc}
\put(545,655){\Carre}
\put(589,655){\Noir}
\put(631,655){\Blanc}

\put(400,690){\line( 1, 0){ 275}}
\put(442,690){\Blanc}
\put(482,690){\Carre}
\put(526,690){\Noir}
\put(568,690){\Blanc}
\put(608,690){\Carre}
\put(652,690){\Noir}

\put(400,725){\line( 1, 0){ 275}}
\put(419,725){\Carre}
\put(463,725){\Noir}
\put(505,725){\Blanc}
\put(545,725){\Carre}
\put(589,725){\Noir}
\put(631,725){\Blanc}

\put(545,520){\makebox(0,0){(c)}}
\end{picture}
\caption{R\'eseaux colori\'es (a) lin\'eaire (b) carr\'e (c) hexagonal}
\label{fig.rescol}
\end{centering}
\end{figure}

Un {\em graphe colori\'e} est un graphe $G$ dont les sommets
sont colori\'es avec $r$ couleurs, not\'ees $1,2,\ldots,r$,
de telle sorte que deux sommets voisins soient de couleurs diff\'erentes;
on notera $G_{s}$ l'ensemble
des sommets de couleur $s\ (1 \leq s \leq r)$;
la \fig{rescol} montre un r\'eseau lin\'eaire
et un r\'eseau carr\'e colori\'es avec deux couleurs,
ainsi qu'un r\'eseau hexagonal colori\'e avec trois couleurs.

Dans un {\em empilement colori\'e} chaque couche $C_i\ (1 \leq i \leq n)$
est un sous-ensemble de $G_{s}$\/, avec $i \equiv s \bmod r$\/;
les couches sont donc successivement de couleurs
$1,2,\ldots r$, et ainsi de suite de fa\c con cyclique.
La d\'efinition d'une coloration garantit que deux
cellules d'une m\^eme couche ne sont
jamais voisines dans $G$.

Dans un empilement standard, chaque cellule qui ne fait pas partie
de la base, repose sur une voisine de la couche imm\'ediatement
inf\'erieure; dans un empilement colori\'e, on impose seulement
qu'une telle cellule poss\`ede au moins une voisine dans l'une
des $r$ couches inf\'erieures. Autrement dit le placement
d'un nouveau sommet $a$ dans un empilement colori\'e
se fait selon le m\^eme principe que dans un empilement standard,
mais en tenant compte de la couleur de $a$ dans le choix de
la couche la plus basse au-dessus des fibres voisines.

%

Bien qu'il n'existe pas de relation simple entre les couches
d'un empilement standard, et celles d'un empilement colori\'e,
les mots obtenus par \'enum\'eration des couches dans un
empilement standard, et dans l'empilement colori\'e correspondant,
sont \'equivalents.

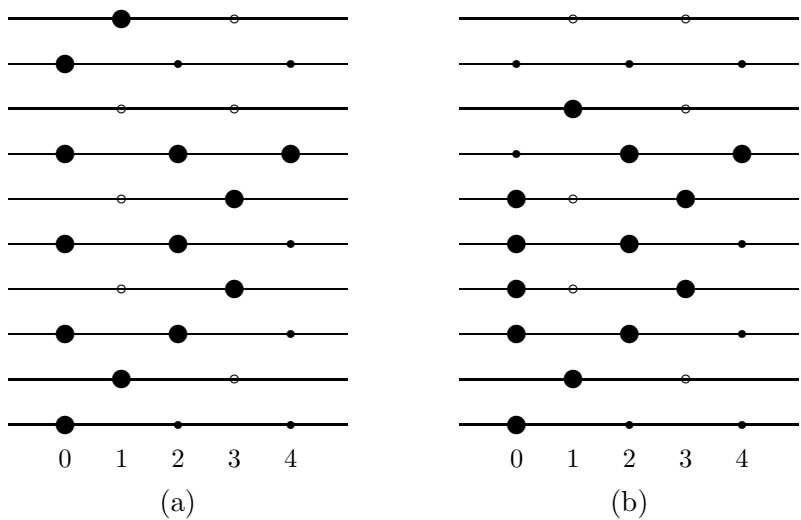
\begin{figure}[htbp]
\begin{centering}
\setlength{\unitlength}{0.015cm}
\newcommand{\Noir}{\circle*{8}}
\newcommand{\Blanc}{\circle{8}}
\newcommand{\Cellule}{\circle*{16}}
\newcommand{\Horizontale}{\line( 1, 0){300}}

\begin{picture}(700,450)(0,-70)
\thinlines


\put(0,0){\Horizontale}
\put(50,0){\Noir}
\put(150,0){\Noir}
\put(250,0){\Noir}

\put(0,40){\Horizontale}
\put(100,40){\Blanc}
\put(200,40){\Blanc}

\put(0,80){\Horizontale}
\put(50,80){\Noir}
\put(150,80){\Noir}
\put(250,80){\Noir}

\put(0,120){\Horizontale}
\put(100,120){\Blanc}
\put(200,120){\Blanc}

\put(0,160){\Horizontale}
\put(50,160){\Noir}
\put(150,160){\Noir}
\put(250,160){\Noir}

\put(0,200){\Horizontale}
\put(100,200){\Blanc}
\put(200,200){\Blanc}

\put(0,240){\Horizontale}
\put(50,240){\Noir}
\put(150,240){\Noir}
\put(250,240){\Noir}

\put(0,280){\Horizontale}
\put(100,280){\Blanc}
\put(200,280){\Blanc}

\put(0,320){\Horizontale}
\put(50,320){\Noir}
\put(150,320){\Noir}
\put(250,320){\Noir}

\put(0,360){\Horizontale}
\put(100,360){\Blanc}
\put(200,360){\Blanc}

\put(50, -30){\makebox(0,0){\small 0}}
\put(100, -30){\makebox(0,0){\small 1}}
\put(150, -30){\makebox(0,0){\small 2}}
\put(200, -30){\makebox(0,0){\small 3}}
\put(250, -30){\makebox(0,0){\small 4}}
\put(150, -70){\makebox(0,0){(a)}}

\put(50,0){\Cellule}
\put(100,40){\Cellule}

\put(50,80){\Cellule}
\put(150,80){\Cellule}
\put(200,120){\Cellule}

\put(50,160){\Cellule}
\put(150,160){\Cellule}
\put(200,200){\Cellule}

\put(50,240){\Cellule}
\put(150,240){\Cellule}
\put(250,240){\Cellule}

\put(50,320){\Cellule}
\put(100,360){\Cellule}


\put(400,0){\Horizontale}
\put(450,0){\Noir}
\put(550,0){\Noir}
\put(650,0){\Noir}

\put(400,40){\Horizontale}
\put(500,40){\Blanc}
\put(600,40){\Blanc}

\put(400,80){\Horizontale}
\put(450,80){\Noir}
\put(550,80){\Noir}
\put(650,80){\Noir}

\put(400,120){\Horizontale}
\put(500,120){\Blanc}
\put(600,120){\Blanc}

\put(400,160){\Horizontale}
\put(450,160){\Noir}
\put(550,160){\Noir}
\put(650,160){\Noir}

\put(400,200){\Horizontale}
\put(500,200){\Blanc}
\put(600,200){\Blanc}

\put(400,240){\Horizontale}
\put(450,240){\Noir}
\put(550,240){\Noir}
\put(650,240){\Noir}

\put(400,280){\Horizontale}
\put(500,280){\Blanc}
\put(600,280){\Blanc}

\put(400,320){\Horizontale}
\put(450,320){\Noir}
\put(550,320){\Noir}
\put(650,320){\Noir}

\put(400,360){\Horizontale}
\put(500,360){\Blanc}
\put(600,360){\Blanc}

\put(450, -30){\makebox(0,0){\small 0}}
\put(500, -30){\makebox(0,0){\small 1}}
\put(550, -30){\makebox(0,0){\small 2}}
\put(600, -30){\makebox(0,0){\small 3}}
\put(650, -30){\makebox(0,0){\small 4}}
\put(550, -70){\makebox(0,0){(b)}}

\put(450,0){\Cellule}
\put(500,40){\Cellule}

\put(450,80){\Cellule}
\put(550,80){\Cellule}

\put(450,120){\Cellule}
\put(600,120){\Cellule}

\put(450,160){\Cellule}
\put(550,160){\Cellule}

\put(450,200){\Cellule}
\put(600,200){\Cellule}

\put(550,240){\Cellule}
\put(650,240){\Cellule}

\put(500,280){\Cellule}

\end{picture}
\caption{Empilements (a) colori\'e (b) standard, sur r\'eseau lin\'eaire}
\label{fig.empcol}
\end{centering}
\end{figure}

La \fig{empcol} montre un empilement colori\'e sur
le r\'eseau lin\'eaire $G$,
et l'empilement standard associ\'e; les sommets de $G$
sont repr\'esent\'es par des entiers: les entiers pairs sont noirs,
les impairs blancs; deux sommets $i$ et $j$ commutent
si $|i-j| \geq 2$. Les mots obtenus par \'enum\'eration des couches
sont respectivement 0102302302401 et 0102030203241,
qui sont bien \'equivalents. La cellule la plus haute de la fibre
$F_1$ repose sur la fibre $F_0$ dans l'empilement colori\'e,
tandis qu'elle repose sur $F_2$ dans l'empilement standard:
le paradoxe n'est qu'apparent, car 0 et 2 commutent.

\subsection{R\'eseaux feuillet\'es}

A un r\'eseau $G$ de dimension $d$,
colori\'e avec $r$ couleurs,
on peut associer un {\em r\'eseau feuillet\'e} $H$, de dimension $d+1$,
d\'efini comme suit:
les sommets de $H$ sont les couples $(a,i)$, pour tout entier $i$
et tout sommet $a$ de $G$,
de couleur $s \equiv i \bmod r$; un arc de $H$ joint $(a,i)$
\`a $(b,j)$ si $a$ et $b$ sont voisins dans $G$, et si
$i < j < i + r$.

Pour une valeur donn\'ee de l'entier $i$, l'ensemble des sommets
$(a,i)$ de H sera appel\'e couche num\'ero i; par d\'efinition,
une couche est monocolore.
Exemples:

\begin{enumerate}
  \item
  Le r\'eseau feuillet\'e associ\'e au r\'eseau lin\'eaire colori\'e
  de la \fig{rescol}-(a), est le r\'eseau carr\'e
  de la \fig{feuilletes}-(a) (les couches sont horizontales).

\begin{figure}[htbp]
\begin{centering}
\setlength{\unitlength}{0.012cm}
\newcommand{\Noir}{\circle*{14}}
\newcommand{\Blanc}{\circle{14}}
\newcommand{\Carre}{
  \makebox(0,0)[r]{\framebox(14,14){}}
}
\newcommand{\Point}{\circle*{2}}

\begin{picture}(940,320)(-20,-40)
\thinlines


\put(40,0){\Noir}
\put(120,0){\Noir}
\put(200,0){\Noir}

\put(0,40){\Blanc}
\put(80,40){\Blanc}
\put(160,40){\Blanc}
\put(240,40){\Blanc}

\put(40,80){\Noir}
\put(120,80){\Noir}
\put(200,80){\Noir}

\put(0,120){\Blanc}
\put(80,120){\Blanc}
\put(160,120){\Blanc}
\put(240,120){\Blanc}

\put(40,160){\Noir}
\put(120,160){\Noir}
\put(200,160){\Noir}

\put(0,200){\Blanc}
\put(80,200){\Blanc}
\put(160,200){\Blanc}
\put(240,200){\Blanc}

\put(20,-20){\line( 1, 1){240}}
\put(100,-20){\line( 1, 1){160}}
\put(180,-20){\line( 1, 1){80}}
\put(-20,20){\line( 1, 1){200}}
\put(-20,100){\line( 1, 1){120}}
\put(-20,180){\line( 1, 1){40}}

\put(220,-20){\line( -1, 1){240}}
\put(140,-20){\line( -1, 1){160}}
\put(60,-20){\line( -1, 1){80}}
\put(260,20){\line( -1, 1){200}}
\put(260,100){\line( -1, 1){120}}
\put(260,180){\line( -1, 1){40}}

\put(120, -40){\makebox(0,0){(a)}}


\put(350,0){\Blanc}
\put(350,150){\Blanc}
\put(500,0){\Blanc}
\put(500,150){\Blanc}

\put(440,30){\Blanc}
\put(440,180){\Blanc}
\put(590,30){\Blanc}
\put(590,180){\Blanc}

\put(470,90){\Noir}

\multiput(350, 0)(10,0){15}{\Point}
\multiput(350, 150)(10,0){15}{\Point}
\multiput(440, 30)(10,0){15}{\Point}
\multiput(440, 180)(10,0){15}{\Point}

\multiput(350, 0)(0,10){15}{\Point}
\multiput(500, 0)(0,10){15}{\Point}
\multiput(440, 30)(0,10){15}{\Point}
\multiput(590, 30)(0,10){15}{\Point}

\multiput(350, 0)(12,4){7}{\Point}
\multiput(350, 150)(12,4){7}{\Point}
\multiput(500, 0)(12,4){7}{\Point}
\multiput(500,150)(12,4){7}{\Point}

\thicklines

\put(350,0){\line(4,3){240}}
\put(350,150){\line( 2, -1){240}}

\put(440,30){\line( 1, 2){60}}
\put(440,180){\line( 1, -3){60}}

\put(470, -40){\makebox(0,0){(b)}}


\put(700,80){\Blanc}
\put(778,160){\Carre}
\put(820,0){\Noir}
\put(900,80){\Blanc}

\put(820,280){\Noir}
\put(738,200){\Carre}
\put(860,120){\Blanc}
\put(938,200){\Carre}

\multiput(700,80)(4,12){10}{\Point}
\multiput(780,160)(4,12){10}{\Point}
\multiput(900,80)(4,12){10}{\Point}

\multiput(700,80)(10,10){8}{\Point}
\multiput(740,200)(10,10){8}{\Point}
\multiput(860,120)(10,10){8}{\Point}

\multiput(740,200)(12,-8){10}{\Point}
\multiput(820,280)(12,-8){10}{\Point}
\multiput(780,160)(12,-8){10}{\Point}

\thicklines

\put(820,0){\line(1, 3){40}}
\put(820,0){\line(3, 5){120}}
\put(820,0){\line(1, 1){80}}

\put(700,80){\line(3, -2){120}}
\put(740,200){\line(2, -5){80}}
\put(780,160){\line(1, -4){40}}

\put(820, -40){\makebox(0,0){(c)}}

\end{picture}
\caption{R\'eseaux feuillet\'es}
\label{fig.feuilletes}
\end{centering}
\end{figure}
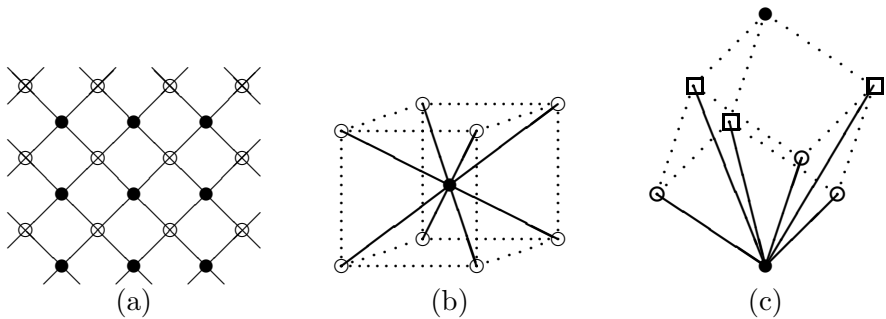

  \item
  Le r\'eseau feuillet\'e associ\'e au r\'eseau carr\'e colori\'e
  de la \fig{rescol}-(b), est le r\'eseau cubique centr\'e;
  on l'obtient en ajoutant aux sommets  d'un r\'eseau cubique,
  colori\'es en blanc, les centres des cubes \'el\'ementaires,
  colori\'es en noir:
  voir la \fig{feuilletes}-(b).
  Ainsi chaque sommet de la couche $i$ poss\'ede 8
  voisins imm\'ediats (les sommets du cube dont il est le centre),
  r\'epartis dans les couches $i-1$ et $i+1$.

  \item
  Le r\'eseau feuillet\'e associ\'e au r\'eseau hexagonal colori\'e
  de la \fig{rescol}-(c), est le r\'eseau cubique {\em nnn} --
  voir \sct{thermo} pour la d\'efinition des r\'eseaux {\em nnn}\/.
  Chaque hexagone du r\'eseau plan est consid\'er\'e comme la
  projection d'un cube ``pos\'e sur la pointe'';
  chaque sommet  de la couche $i$,
  poss\`ede 3 voisins imm\'ediats dans la couche $i+1$,
  et 3 ``voisins suivants'' dans la couche $i+2$:
  voir la \fig{feuilletes}-(c).

\end{enumerate}

Un r\'eseau feuillet\'e peut \^etre {\em orient\'e}
en dirigeant les arcs des couches inf\'erieures vers
les couches sup\'erieures (il n'y a jamais d'ar\^ete
\`a l'int\'erieur d'une couche); il peut \^etre
{\em \'etendu} en ajoutant les arcs
$(a,i) \rightarrow (a,i+r)$.
Par exemple le r\'eseau feuillet\'e \'etendu
associ\'e au r\'eseau lin\'eaire est un r\'eseau triangulaire,
obtenu en ajoutant des arcs verticaux \`a la \fig{feuilletes}-(a).

\subsection{Equivalence entre animaux dirig\'es et empilements}
\label{sec.equivalence}

\begin{theoreme}
  \label{th.equivalence}
  Soit $G$ un r\'eseau colori\'e, de dimension $d$,
  $H$ le r\'eseau feuillet\'e orient\'e issu de $G$,
  de dimension $d+1$, et $H'$ le r\'eseau \'etendu.
  Il existe une bijection entre pyramides
  (resp. pyramides strictes) construites sur $G$,
  et animaux dirig\'es, de sources ponctuelles,
  sur $H'$ (resp. $H$).
\end{theoreme}

\preuve soit $r$ le nombre de couleurs de $G$;
par d\'efinition d'un empilement colori\'e $E$,
chaque cellule $(b,j)$ qui n'appartient pas \`a la base
de l'empilement, poss\`ede une voisine $(a,i)$ dans
l'une des $r$ couches inf\'erieures.
Si $j=i+r$, alors $a=b$, car les couches $i$ et $i+r$
ont m\^eme couleur; donc si l'empilement est strict,
on a $j < i+r$. Par d\'efinition d'un r\'eseau
feuillet\'e orient\'e, il existe un arc
$(a,i) \rightarrow (b,j)$ dans $H'$, et m\^eme dans $H$
si l'empilement est strict.

De proche en proche, on en d\'eduit que, pour toute cellule de
$E$, il existe un chemin dans $H$ (ou dans $H'$) depuis
la base de $E$ jusqu'\`a cette cellule.
Un empilement colori\'e sur $G$ est donc un animal dirig\'e
sur $H'$, et m\^eme sur $H$ s'il s'agit d'un empilement
strict; et vice-versa. Si la source de l'animal
est ponctuelle, la base de l'empilement est r\'eduite \`a
un point: c'est donc une pyramide.

Enfin on a vu \sct{empilements-colories} qu'il existe une
bijection naturelle entre empilements standard
et empilements colori\'es\cqfd
\bigskip

\pagebreak[2]

Exemples:
\begin{enumerate}
  \item
  Lorsque le nombre de couleurs $r$ vaut 2, il n'existe aucune
  diff\'erence entre empilements stricts standard
  et empilements stricts colori\'es. Donc l'animal
  dirig\'e de la \fig{animaux}-(a) est exactement une pyramide
  stricte sur r\'eseau lin\'eaire: il suffit de faire subir
  au dessin une rotation de $\pi/4$ dans le sens trigonom\'etrique
  pour le v\'erifier.
  
  \item
  L'animal  dirig\'e de la \fig{animaux}-(b) est une pyramide
  colori\'ee (non stricte) sur r\'eseau lin\'eaire: il faut
  faire subir cette fois
  au dessin une rotation de $\pi/6$ dans le sens trigonom\'etrique.
  La pyramide colori\'ee de la \fig{empcol}-(a)
  est un autre exemple d'animal dirig\'e sur r\'eseau triangulaire;
  la pyramide standard de la \fig{empcol}-(b)
  est donc la pyramide associ\'ee \`a cet animal.
  
  \item
  Une pyramide stricte sur r\'eseau carr\'e correspond \`a un animal
  dirig\'e sur un r\'eseau cubique centr\'e;
  l'animal peut se d\'evelopper,
  dans un rep\`ere convenablement choisi,
  selon l'un des 4 vecteurs:
  \[
    (1,1,1),\ (1,-1,1),\ (-1,1,1),\ (-1,-1,1)
  \]
  auxquels il faut ajouter le vecteur $(0,0,2)$ pour traiter
  les pyramides quelconques; voir \fig{feuilletes}-(b).
  
  \item
  Une pyramide stricte sur r\'eseau hexagonal correspond \`a un animal
  dirig\'e sur un r\'eseau cubique {\em nnn}\/;
  l'animal peut se d\'evelopper selon l'un des 6 vecteurs:
  \[
    (1,0,0),\ (0,1,0),\ (0,0,1),\
    (1,1,0),\ (1,0,1),\ (0,1,1)
  \]
  auxquels il faut ajouter le vecteur $(1,1,1)$ pour traiter
  les pyramides quelconques.
  Noter qu'avec le rep\`ere choisi ici, les couches ne sont
  plus horizontales, comme sur la \fig{feuilletes}-(c),
  mais d'\'equations $x+y+z=i$.

\end{enumerate}

En combinant ce th\'eor\`eme, d\^u \`a \viennot,
avec le \cor{densite}
de la \sct{applications-gaz}, et avec le \thref{stricts},
qui relie les s\'eries caract\'eristiques des empilements 
selon qu'ils sont stricts ou non, 
on retrouve le r\'esultat suivant de
D.Dhar \cite{dhar2}:

\begin{corollaire}
  Le probl\`eme de l'\'enum\'eration des animaux dirig\'es
  sur un r\'eseau feuillet\'e $H$ issu d'un r\'eseau $G$,
  est \'equivalent \`a celui du
  calcul de la densit\'e moyenne d'un gaz \`a particules dures sur $G$,
  en fonction de la temp\'erature\cqfd
\end{corollaire}

Il est int\'eressant de comparer les m\'ethodes employ\'ees
en physique th\'eo\-rique et celles, combinatoires, pr\'esent\'ees
dans cet article. Dans \cite{dhar2}, D.Dhar emploie l'approche
suivante:
\begin{quote}
  ``
    Un site $j$ est appel\'e successeur du site $i$
    s'il existe un arc de $i$ vers $j$.
    Pour un ensemble de sites $C$, l'ensemble des successeurs
    sera not\'e $S(C)$.
     (\ldots)
    D\'efinissons pour chaque site $i$ du r\'eseau
    un temps $t(i)$, \`a valeurs enti\`eres,
    tel que $t(j) \geq t(i)+1$ chaque fois que $j$
    est un successeur de $i$.
    (Dans un animal dirig\'e)
    les configurations autoris\'ees de sites occup\'es,
    sur la surface $t=T$, ne d\'ependent que des configurations
    aux instants pr\'ec\'edents $t \leq T-1$.
    Consid\'erons d'abord, pour simplifier, le cas o\`u elles
    d\'ependent seulement de la configuration sur la surface
    $t=T-1$ (cas du r\'eseau carr\'e, et du r\'eseau cubique centr\'e).
    Soit $C$ l'ensemble des sites occup\'es sur une hyper-surface
    $t=c^{\rm ste}$. D\'efinissons les fonctions g\'en\'eratrices
    $A_C(x)$, qui \'enum\`erent selon leur taille les
    animaux dirig\'es de source $C$.
    
    La propri\'et\'e markovienne des animaux dirig\'es implique
    que ces fonctions g\'en\'eratrices v\'erifient les r\'ecurrences:
    \[
    A_C (x) = x^{|C|}\left(1+\sum_D A_D(x) \right)
    \]
    o\`u la sommation est \'etendue \`a tous les sous-ensembles
    non vides $D$ de $S(C)$. Notons que le temps sur les sites de $D$
    est sup\'erieur de 1 au temps sur les sites de $C$.
    Il est facile de g\'en\'eraliser
    ces r\'ecurrences au cas o\`u les configurations autoris\'ees
    \`a l'instant $t=T$ d\'ependent aussi des configurations
    aux instants $t \leq T-2$ (cas des r\'eseaux {\em nnn}, etc\ldots).
    Ces \'equations forment la hi\'erarchie de
    Bogoliubov-Born-Green-Kirkwood-Yvon (sic) pour le probl\`eme
    de l'\'enum\'eration des animaux dirig\'es.
  ''
\end{quote}

D.Dhar compare alors ces \'equations \`a celles qu'on obtient
en simulant l'\'evolution d'un gaz: \`a chaque incr\'ementation
du temps, chaque site \'evolue al\'eatoirement,
la \proba\ d'occupation d'un site donn\'e ne d\'ependant que de la
configuration des sites voisins, avec une formule simple
d\'eduite de la loi de Gibbs \eq{proba}; lorsque
$t \rightarrow \infty$, le gaz approche ainsi de l'\'equilibre
thermique. Ce principe de simulation, classique en physique
statistique, est aussi \`a la base des r\'eseaux de neurones
(voir par exemple \cite{weisbuch}).
\begin{quote}
  ``
    (Etudions) l'\'evolution stochastique dans le temps d'un gaz
    sur un r\'eseau de dimension $d-1$. Ce r\'eseau est constitu\'e
    de $r$ sous-r\'eseaux entrelac\'es $L_1, L_2, \ldots,L_r$,
    chacun d'entre eux \'etant isomorphe \`a un hyperplan
    $\tau=c^{\rm ste}$ du r\'eseau de dimension $d$ (sur lequel
    se d\'eveloppent les animaux dirig\'es). (\ldots)
    A l'instant $\tau=mr+i \; (1\leq i \leq r)$,
    la configuration des sites occup\'es sur le sous-r\'eseau
    $L_i$ \'evolue al\'eatoirement, tandis que les autres
    sous-r\'eseaux restent inchang\'es. Lorsque $\tau \rightarrow \infty$,
    la \proba\ des diff\'erentes configurations du gaz
    tend vers une distribution limite invariante. Pour des lois
    de transition arbitraires, le calcul de la distribution
    invariante est tout \`a fait difficile. (Ici) les lois de
    transition sont particuli\`erement simples, et correspondent
    aux algorithmes de Monte-Carlo utilis\'es couramment
    pour \'etudier les lois d'\'equilibre d'un gaz sur un r\'eseau
    avec particules dures (\ldots)
    La \proba\ d'occupation d'un site $i$ vaut 0 si 
    l'un des sites successeurs de $i$ est occup\'e;
    et elle vaut $p$ s'ils sont tous vides.
  ''
\end{quote}
D.Dhar remarque qu'il suffit de ``renverser le temps'',
c'est \`a dire de poser $\tau=-t$, pour faire dispara\^\i tre
le paradoxe qui fait d\'ependre la \proba\ d'occupation d'un site $i$
de la configuration de ses successeurs,
au lieu de ses pr\'ed\'ecesseurs. Il \'ecrit:
\begin{quote}
  ``
    Soit $P(C)$ la \proba\ que tous les sites de 
    l'ensemble $C$, inclus dans une hyper-surface
    $\tau=c^{\rm ste}$,  soient occup\'es.
    Je suppose $\tau$ suffisamment grand pour que $P(C)$
    ne d\'epende pas de $\tau$.
    Cette \proba\ est \'egale \`a $p^{|C|}$ multipli\'e par
    la \proba\ que tous les sites de $S(C)$ soient vides.
    En utilisant le principe d'inclusion-exclusion
    \[
      P(C) = p^{|C|}\left(1+\sum_D(-1)^{|D|} P(D) \right)
    \]
    o\`u la sommation est \'etendue \`a tous les sous-ensembles
    non vides $D$ de $S(C)$. On voit donc que
    \[
      A_C (x=-p) = (-1)^{|C|}P(C).
    \]
    Ainsi la fonction g\'en\'eratrice des animaux dirig\'es
    issus d'un seul point source
    (sur un r\'eseau de dimension $d$)
    est d\'etermin\'ee par la
    densit\'e moyenne de sites occup\'es dans un gaz \`a
    l'\'equilibre sur un r\'eseau de dimension $d-1$.
  ''
\end{quote}

On voit comment le \theor\ d'inversion est remplac\'e par
l'utilisation du principe d'inclusion-exclusion
en calcul des \proba s. L'introduction de sous-r\'eseaux
entrelac\'es, qui a \'et\'e utile aussi dans
le pr\'esent article pour \'eclaicir la bijection entre
empilements et animaux, n'est pas expliqu\'ee par D.Dhar.
Elle est rendue n\'ecessaire par le mod\`ele
des particules dures,
car si \`a l'instant $\tau$ tous
les sites changent al\'eatoirement d'\'etat,
on peut voir
appara\^\i tre des particules sur des sites voisins,
chacun \'etant
pourtant entour\'e de sites vides \`a l'instant pr\'ec\'edent.

Nos m\'ethodes alg\'ebriques sont donc parentes
de celles des physiciens, dont on ne peut qu'admirer
l'ing\'eniosit\'e; mais l'alg\`ebre et la combinatoire,
substitu\'ees aux r\'ecurrences d\'eduites du
calcul des probabilit\'es,
permettent d'\'eliminer l'impr\'ecision de certains arguments;
surtout elles
permettent une \'etude plus fine des objets: dans la section
suivante, l'alg\`ebre des empilements nous donnera
des d\'ecompositions canoniques des
animaux en dimension 2, qui permettent \`a la fois
de retrouver  tr\`es simplement les s\'eries g\'en\'eratrices
qui les \'enum\`erent, et de construire des animaux al\'eatoires,
avec distribution uniforme, en temps lin\'eaire.

\newcommand{\binome}[2]{
  \left( \hspace{-0.3em}
    \begin{array}{c}
      #1 \\ #2
    \end{array}
  \hspace{-0.3em} \right)
}

\newcommand{\Mz}{{\bf M}}

\section{Solution du mod\`ele lin\'eaire}

\label{sec.solution}

\subsection{D\'ecomposition en \eqr s}

Le cas le plus simple de mod\`ele de gaz \`a particules dures
sur un r\'eseau $G$, est celui o\`u $G$ est de dimension 1,
autrement dit o\`u $G$ est une cha\^\i ne infinie,
isomorphe \`a {\bf Z}.
D'apr\`es la \sct{equivalence}, les pyramides strictes
sur $G$ sont exactement les animaux dirig\'es, de source ponctuelle,
sur r\'eseau carr\'e.

\begin{figure}[htbp]
\begin{centering}
\setlength{\unitlength}{0.015cm}
\begin{picture}(420,400)(-40,-20)
\thicklines


\put(0,100){\line(1,-2){50}}
\put(0,100){\line(1, 0){50}}
\put(50,0){\line(0,1){100}}
\put(35,70){\makebox(0,0){L}}

\thinlines
\put(50,100){\line(1,1){30}}
\thicklines

\put(-20,330){\line(1,-2){100}}
\put(-20,330){\line(1, 0){200}}
\put(80,130){\line(1,2){100}}
\put(80,260){\makebox(0,0){Q}}

\multiput(-40,-20)(10,0){4}{\circle*{2}}
\put(20,-20){\makebox(0,0){\small -1}}
\put(50,-20){\makebox(0,0){\small 0}}
\put(80,-20){\makebox(0,0){\small 1}}
\put(110,-20){\makebox(0,0){\small 2}}
\multiput(140,-20)(10,0){4}{\circle*{2}}


\put(250,100){\line(1,-2){50}}
\put(250,100){\line(1, 0){50}}
\put(300,0){\line(0,1){100}}

\thinlines
\put(300,100){\line(1,1){30}}
\thicklines

\put(280,230){\line(1,-2){50}}
\put(280,230){\line(1, 0){50}}
\put(330,130){\line(0,1){100}}

\multiput(330,230)(10,10){5}{\circle*{2}}

\put(330,380){\line(1,-2){50}}
\put(330,380){\line(1, 0){50}}
\put(380,280){\line(0,1){100}}

\multiput(210,-20)(10,0){4}{\circle*{2}}
\put(270,-20){\makebox(0,0){\small -1}}
\put(300,-20){\makebox(0,0){\small 0}}
\put(330,-20){\makebox(0,0){\small 1}}
\multiput(350,-20)(10,0){2}{\circle*{2}}
\put(380,-20){\makebox(0,0){\small i}}
\multiput(410,-20)(10,0){4}{\circle*{2}}

\end{picture}
\caption{Sch\'ema de \decomp\ d'une pyramide en produit d'\eqr s}
\label{fig.schema-pyreq}
\end{centering}
\end{figure}
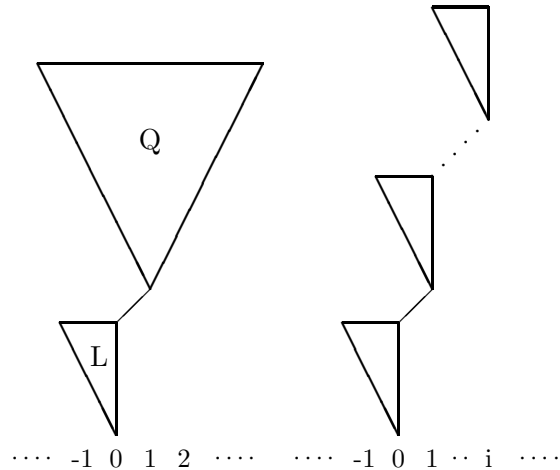

Rep\'erons les sommets de $G$ par des entiers,
et consid\'erons une pyramide $P$ de base~0.
On peut d\'ecomposer $P$ en remarquant que:
\begin{itemize}
  \item
  si $P$ ne contient aucune cellule au-dessus de 1,
  autrement dit si la fibre $F_1$ est vide,
  alors toutes les fibres $F_i$, pour $i>1$,
  sont vides aussi; on dira dans ce cas que $P$
  est une {\em\eqr} de base~0.
  Cette propri\'et\'e repose sur une particularit\'e
  \'evidente du mod\`ele lin\'eaire: il n'existe pas
  de chemin de 0 vers $i>1$, qui ne passe pas par~1.
  
  \item
  si $F_1$ n'est pas vide, soit $(1,k)$ la cellule la plus
  basse sur cette fibre, et soit $Q$ la pyramide engendr\'ee
  par cette cellule; la \decomp\ \eq{decomposition}
  de la \sct{pyramides} devient:
  \begin{equation}
    \label{eq.pyreq}
    P = L\: Q
  \end{equation}
  o\`u L est une \'equerre, car $Q$ contient la fibre $F_1$.
\end{itemize}
La \decomp\ \eq{pyreq} peut \'evidemment
\^etre it\'er\'ee, en d\'ecomposant $Q$, etc\ldots;
elle est sch\'ematis\'ee, ainsi que son it\'eration,
sur la \fig{schema-pyreq}.
La plus grande valeur de $i$ pour laquelle la fibre $F_i$
n'est pas vide, est appel\'ee
{\em \dld} de la pyramide (ou de l'animal) $P$;
le nombre d'\eqr s qui interviennent
dans la d\'ecomposition de $P$ est \'egal \`a $i+1$;
une \eqr\ est une pyramide de \dld\ nulle.

\begin{figure}[htbp]
\begin{centering}
\setlength{\unitlength}{0.02cm}
\newcommand{\Noir}{\circle*{10}}
\newcommand{\Blanc}{\circle{10}}
\begin{picture}(580,300)(0,-20)
\thinlines


\multiput(140,0)(-20,20){8}{\line(1,1){140}}
\multiput(0,140)( 20,20){8}{\line(1,-1){140}}

\multiput(140,0)(-20,20){6}{\Noir}
\multiput(140,40)(-40,40){3}{\Noir}
\multiput(120,100)(-20,20){4}{\Noir}
\put(80,180){\Noir}
\put(100,200){\Noir}

\put(127,120){\makebox(0,0){\small 1}}
\put(107,220){\makebox(0,0){\small 2}}
\put(127,240){\makebox(0,0){\small 3}}
\put(147,260){\makebox(0,0){\small 4}}

\multiput(160,60)(20,20){3}{\Blanc}
\multiput(180,120)(20,20){3}{\Blanc}
\multiput(140,120)(20,20){4}{\Blanc}
\multiput(180,200)(20,20){2}{\Blanc}
\multiput(120,220)(20,20){3}{\Blanc}

\put(140,-20){\makebox(0,0){(a)}}


\multiput(440,0)(-20,20){8}{\line(1,1){140}}
\multiput(300,140)(20,20){8}{\line(1,-1){140}}

\multiput(460,60)(-20,20){3}{\Blanc}
\put(480,80){\Noir}
\multiput(500,100)(-20,20){4}{\Blanc}
\multiput(500,140)(-20,20){2}{\Blanc}
\multiput(520,160)(-20,20){4}{\Noir}
\put(500,220){\Noir}

\put(427,80){\makebox(0,0){\small 1}}
\put(407,100){\makebox(0,0){\small 2}}
\put(427,160){\makebox(0,0){\small 3}}
\put(447,220){\makebox(0,0){\small 4}}

\put(440,-20){\makebox(0,0){(b)}}

\end{picture}
\caption{Exemple de \decomp\ d'une pyramide en produit d'\eqr s}
\label{fig.pyreq}
\end{centering}
\end{figure}
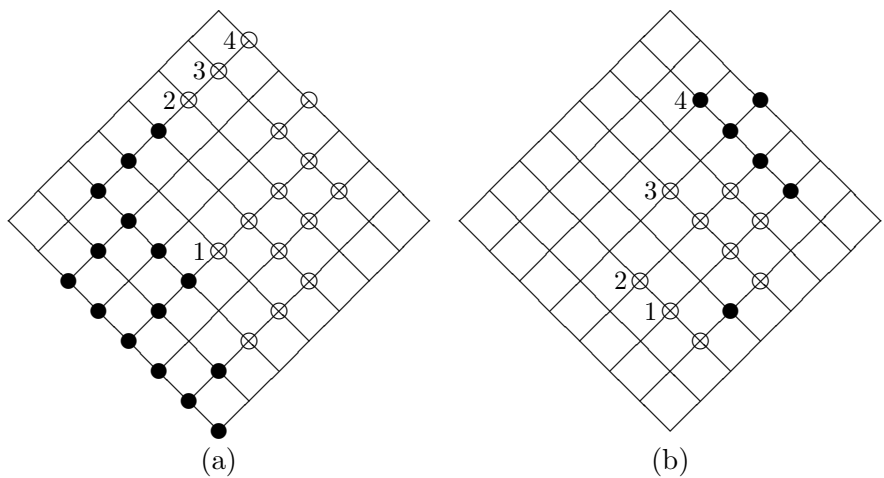

La \fig{pyreq} illustre ces \decomp s sur
l'exemple de l'animal de la \fig{animaux}-(a),
dont la \dld\ vaut 4;
la partie (a) de la figure montre la premi\`ere
\decomp, avec les cellules de $L$ marqu\'ees en noir,
et celles de $Q$ en blanc;
la partie (b) montre  la \decomp\ it\'er\'ee
de $Q$ comme produit de quatre \'equerres.
 On notera
que la pyramide $Q$ ne semble pas former un animal sur la
\fig{pyreq}-(a); (b)
montre $Q$ sous forme d'animal, apr\`es la ``chute''
de quatre cellules, une fois l'\eqr\ $L$ disparue.
Cet exemple
montre combien l'alg\`ebre des empilements est
essentielle pour ces \decomp s,
invisibles si l'on s'en tient \`a la d\'efinition
brute d'un animal.

En notant $\Pi_i$ (resp. $\Lambda_i$) la s\'erie caract\'eristique
des pyramides (resp. \eqr s) de base~$i$, la \decomp\
\eq{pyreq} devient (en tenant compte du cas o\`u $P$
est d\'ej\`a une \'equerre):
  \begin{equation}
    \label{eq.car-pyreq}
    \Pi_0 = \Lambda_0 \: (1 + \Pi_1)
  \end{equation}

Une \eqr\ $L$ de taille strictement sup\'erieure \`a 1,
peut \^etre \`a son tour d\'ecompos\'ee;
soit $x$ la base de $L$:

\begin{itemize}
  \item
  si $x$ est seule sur la fibre $F_0$,
  alors la pyramide $M$, engendr\'ee
  par la cellule $(-1,1)$, est une \eqr, et:
  \begin{equation}
  \label{eq.eqr1}
     L = x\: M
  \end{equation}
  
  \item
  sinon, soit $(0,k)$ la cellule la plus basse sur $F_0$,
  avec $k>1$; cette cellule engendre une pyramide $N$,
  qui est \`a son tour une \eqr, et la \decomp:
  \begin{equation}
  \label{eq.eqr2}
     L = x\: M \: N
  \end{equation}
  laisse une pyramide $M$, de base $(-1,1)$,
  dont on a vid\'e la fibre $F_0$\/; $M$ est donc une \eqr.
\end{itemize}

\begin{figure}[htbp]
\begin{centering}
\setlength{\unitlength}{0.018cm}
\newcommand{\Noir}{\circle*{10}}
\newcommand{\Blanc}{\circle{10}}
\begin{picture}(600,320)(-10,-50)


\thinlines
\put(50,30){\line(1,-1){30}}
\thicklines

\put(0,130){\line(1,-2){50}}
\put(0,130){\line(1, 0){50}}
\put(50,30){\line(0,1){100}}
\put(80,0){\Noir}
\put(80,15){\makebox(0,0){$x$}}
\put(35,100){\makebox(0,0){M}}

\multiput(-10,-20)(10,0){4}{\circle*{2}}
\put(50,-20){\makebox(0,0){\small -1}}
\put(80,-20){\makebox(0,0){\small 0}}
\put(110,-20){\makebox(0,0){\small 1}}
\multiput(140,-20)(10,0){4}{\circle*{2}}

\put(80,-50){\makebox(0,0){(a)}}


\thinlines
\put(250,30){\line(1,-1){30}}
\thicklines

\put(200,130){\line(1,-2){50}}
\put(200,130){\line(1, 0){50}}
\put(250,30){\line(0,1){100}}
\put(280,0){\Noir}
\put(280,15){\makebox(0,0){$x$}}
\put(235,100){\makebox(0,0){M}}

\thinlines
\put(250,130){\line(1,1){30}}
\thicklines

\put(230,260){\line(1,-2){50}}
\put(230,260){\line(1, 0){50}}
\put(280,160){\line(0,1){100}}
\put(265,230){\makebox(0,0){N}}

\multiput(190,-20)(10,0){4}{\circle*{2}}
\put(250,-20){\makebox(0,0){\small -1}}
\put(280,-20){\makebox(0,0){\small 0}}
\put(310,-20){\makebox(0,0){\small 1}}
\multiput(340,-20)(10,0){4}{\circle*{2}}

\put(280,-50){\makebox(0,0){(b)}}


\thinlines

\multiput(500,0)(-20,20){8}{\line(1,1){100}}
\multiput(360,140)( 20,20){6}{\line(1,-1){140}}

\multiput(500,0)(-20,20){6}{\Noir}
\multiput(460,80)(-40,40){2}{\Noir}
\multiput(480,100)(-20,20){4}{\Blanc}
\put(500,40){\Blanc}
\put(440,180){\Blanc}
\put(460,200){\Blanc}

\put(500,-50){\makebox(0,0){(c)}}

\end{picture}
\caption{D\'ecomposition d'une \eqr\ en produit d'\eqr s}
\label{fig.equerres}
\end{centering}
\end{figure}
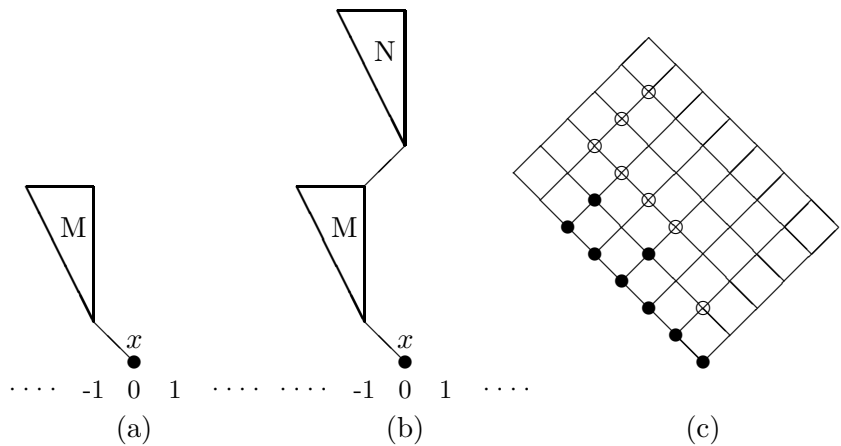

Ces \decomp s sont sch\'ematis\'ees
sur la \fig{equerres}, parties (a) et (b);
la partie (c) illustre (b) sur l'exemple
de l'\eqr\ noire de la \fig{pyreq}-(a);
l'\eqr\ $M$ est laiss\'ee noire,
et les cellules de $N$ apparaissent en blanc.
On peut voir \`a nouveau que cette \decomp\
n'est pas visible sans consid\'erer un animal comme un empilement,
car $N$ n'est pas connexe.

Les \decomp s \eq{eqr1} et \eq{eqr2}
entra\^\i nent, pour les s\'eries caract\'eristi\-ques,
en tenant compte du cas o\`u $L=x$ (pyramide r\'eduite \`a sa base):
  \begin{equation}
    \label{eq.car-eqr-strict}
    \Lambda_0 = x (1+ \Lambda_{-1} + \Lambda_{-1} \: \Lambda_0)
  \end{equation}

Les \'equations \eq{car-pyreq} et \eq{car-eqr-strict} concernent
les pyramides strictes, ou les animaux dirig\'es sur r\'eseau
carr\'e; pour les pyramides g\'en\'erales, qui correspondent
aux animaux dirig\'es sur r\'eseau triangulaire,
on peut appliquer le \thref{stricts},
ou remarquer que \eq{car-pyreq} reste valable, et que
\eq{car-eqr-strict} devient:
  \begin{equation}
    \label{eq.car-eqr}
    \Lambda_0 = x (1+ \Lambda_{-1} + \Lambda_0 + \Lambda_{-1} \: \Lambda_0)
  \end{equation}
car il y a une seconde \decomp\ de type \ref{fig.equerres}-(a),
avec $M$, de base $0$, reposant verticalement sur $x$.

\subsection{S\'eries g\'en\'eratrices}

On obtient les s\'eries g\'en\'eratrices des pyramides et des
\eqr s en projetant (voir \sct{applications-gaz}) les \'equations
obtenues pour les s\'eries caract\'eristiques dans la section
pr\'ec\'edente; ces projections ne d\'ependent pas de la base.
En notant $\pi(t)$ (resp. $\lambda(t)$)
la \sgen\ qui \'enum\`ere les
pyramides strictes (resp. les \eqr s strictes),
\eq{car-pyreq} et \eq{car-eqr-strict} fournissent:
  \begin{eqnarray}
     \label{eq.gen-pyreq}
     \pi(t) & = & \lambda(t) \: [1 + \pi(t)] \\
     \label{eq.gen-eqr-strict}
     \lambda(t) & = & t \: [1 + \lambda(t) + \lambda(t) ^2 ]
  \end{eqnarray}

Ces \'equations sont facilement r\'esolues:
  \begin{eqnarray}
    \label{eq.sol-eqr-strict}
    \lambda(t) & = & \frac{1-t-\sqrt{(1+t)(1-3t)}}{2t} \\
    & = & t + t^2  + 2 t^3  + 4 t^4  + 9 t^5 
            + 21 t^6  + 51 t^7  + 127 t^8  + 323 t^9 
            + \ldots \nonumber \\
    \nonumber \\
    \label{eq.sol-pyr-strict}
    \pi(t) & = & \frac{1}{2}
    \left( \sqrt{\frac{1+t}{1-3t}} - 1\right) \\
    & = &  t + 2 t^2  + 5 t^3  + 13 t^4  + 35 t^5
             + 96 t^6  + 267 t^7  + 750 t^8  + 2123 t^9 
             + \ldots \nonumber
  \end{eqnarray}
  
Pour \'enum\'erer les \eqr s et les pyramides g\'en\'erales,
il faut remplacer \eq{gen-eqr-strict} par l'\'equation
d\'eduite de \eq{car-eqr}:
  \begin{equation}
     \label{eq.gen-eqr}
     \lambda(t)  =  t \: [1 + 2 \lambda(t) + \lambda(t) ^2]
  \end{equation}
et on obtient cette fois:
  \begin{eqnarray}
    \label{eq.sol-eqr}
    \lambda(t) & = & \frac{1-2t-\sqrt{1-4t}}{2t}
    \;\; = \;\; \sum \frac{1}{n+1} \binome{2n}{n} t^n \\
    & = & t + 2 t^2  + 5 t^3  + 14 t^4  + 42 t^5 
            + 132 t^6  + 429 t^7  + 1430 t^8 
            + \ldots \nonumber \\
    \nonumber \\
    \label{eq.sol-pyr}
    \pi(t) & = & \frac{1}{2}
    \left( \frac{1}{\sqrt{1-4t}} - 1\right)
    \;\; = \;\; \frac{1}{2}\sum \binome{2n}{n} t^n \\
    & = &  t + 3 t^2  + 10 t^3  + 35 t^4  + 126 t^5
             + 462 t^6  + 1716 t^7  + 6435 t^8 
             + \ldots \nonumber
  \end{eqnarray}
  
Il n'existe pas de forme close pour les coefficients des s\'eries
\eq{sol-eqr-strict} et \eq{sol-pyr-strict},
contrairement au cas des s\'eries
\eq{sol-eqr} et \eq{sol-pyr}.
On v\'erifie que, conform\'ement au \thref{stricts}, on peut passer
du premier couple de s\'eries au second
par la substitution
$t \rightarrow t/(1-t)$, et inversement
du second au premier par la substitution
$t \rightarrow t/(1+t)$.

En appliquant le \thref{equivalence}, \sct{equivalence},
on retrouve les r\'esultats de \cite{dhar2,gouyou,hakim,nadal}:

\begin{theoreme}
  \label{th.enum-animaux}
  La s\'erie \eq{sol-pyr-strict} -- resp. \eq{sol-pyr} -- est la \senum\
  des animaux dirig\'es, de source ponctuelle, sur r\'eseau
  carr\'e --~resp. triangulaire~--~, selon leur taille\cqfd
\end{theoreme}

Si, d'autre part, on applique le \cor{densite}
de la \sct{applications-gaz}, on obtient:
\begin{theoreme}
  \label{th.densite-lineaire}
  La densit\'e moyenne d'un gaz avec particules dures,
  sur r\'e\-seau lin\'eaire, vaut:
  \[
    d(t) = \frac{1}{2}\left( 1 - \frac{1}{\sqrt{1+4t}} \right)
  \]  
\end{theoreme}

\preuve d'apr\`es le \cor{densite},
cette densit\'e vaut $- \pi(-t)$, o\`u $\pi(t)$
est la \senum\ des pyramides, donn\'ee par \eq{sol-pyr}\cqfd
\bigskip

Rappelons que la variable $t$ est li\'ee \`a la temp\'erature absolue $T$
par $t=e^{\mu/kT}$, o\`u $\mu$ d\'esigne le potentiel
chimique, et $k$ la constante de Boltzmann (cf. \sct{particules-dures}).
Si le potentiel chimique est positif, l'\'energie diminue lorsque
le nombre de particules augmente;
lorsque la temp\'erature $T$ cro\^\i t de $0$ \`a $+ \infty$,
$t$ d\'ecro\^\i t de $+ \infty$ \`a $1$, et
$d(t)$ d\'ecro\^\i t de $1/2$
(configuration d'\'energie minimale:
un site sur deux est occup\'e) 
\`a $(1-1/\sqrt{5})/2$
(toutes les configurations sont \'equiprobables).
On note que la densit\'e ne
pr\'esente pas de singularit\'e, bien que son d\'eveloppement en s\'erie
ait pour rayon de convergence $1/4$.

D.Dhar, dans \cite{dhar1,dhar2}, d\'eduit le \thref{enum-animaux}
du \thref{densite-lineaire}, bien connu des physiciens, alors
que notre approche donne les deux r\'esultats simultan\'ement.

\subsection{Chemins de Dyck et chemins de Motzkin}
\label{sec.motzkin}

\begin{figure}[htbp]
\begin{centering}
\setlength{\unitlength}{0.018cm}
\newcommand{\Montant}{\line(1,1){20}}
\newcommand{\Descendant}{\line(1,-1){20}}
\newcommand{\Horizontal}{\line(1,0){20}}
\begin{picture}(600,140)(0,-30)


\thinlines
\put(0,0){\vector(1,0){310}}
\put(305,15){\makebox(0,0){$x$}}
\multiput(0,20)(0,20){4}{\line(1,0){280}}
\put(0,0){\vector(0,1){110}}
\put(15,105){\makebox(0,0){$y$}}
\multiput(20,0)(20,0){14}{\line(0,1){80}}
\put(140,-30){\makebox(0,0){(a)}}

\thicklines
\multiput(0,0)(20,20){2}{\Montant}
\put(40,40){\Descendant}
\put(60,20){\Montant}
\multiput(80,40)(20,-20){2}{\Descendant}
\put(120,0){\Montant}
\put(140,20){\Descendant}
\multiput(160,0)(20,20){3}{\Montant}
\multiput(220,60)(20,-20){3}{\Descendant}


\thinlines
\put(340,0){\vector(1,0){310}}
\put(645,15){\makebox(0,0){$x$}}
\multiput(340,20)(0,20){4}{\line(1,0){280}}
\put(340,0){\vector(0,1){110}}
\put(355,105){\makebox(0,0){$y$}}
\multiput(360,0)(20,0){14}{\line(0,1){80}}
\put(480,-30){\makebox(0,0){(b)}}

\thicklines
\put(340,0){\Montant}
\put(360,20){\Horizontal}
\multiput(380,20)(20,20){2}{\Montant}
\multiput(420,60)(20,-20){3}{\Descendant}
\multiput(480,0)(20,0){2}{\Horizontal}
\put(520,0){\Montant}
\put(540,20){\Descendant}
\put(560,0){\Montant}
\put(580,20){\Horizontal}
\put(600,20){\Descendant}

\end{picture}
\caption{Chemins (a) de Dyck (b) de Motzkin}
\label{fig.ex-chemins}
\end{centering}
\end{figure}
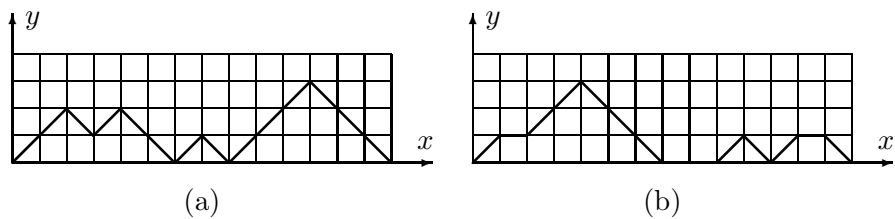

Les formules de la section pr\'ec\'edente sont classiques
en combinatoire, dans d'autres contextes. Rappelons qu'un
{\em chemin de Dyck} est un chemin dans le plan discret,
compos\'e de pas ascendants $(1,1)$ et de pas descendants
$(1,-1)$, et qui joint l'origine $(0,0)$ \`a un point
de l'axe des abscisses, sans quitter le quart de plan
sup\'erieur: voir \fig{ex-chemins}-(a).
Un {\em chemin de Motzkin} est d\'efini de fa\c con semblable,
mais comporte en outre des pas horizontaux $(1,0)$:
voir \fig{ex-chemins}-(b).

Dans un chemin de Motzkin {\em r-colori\'e},
il y a $r$ sortes de pas horizontaux;
dans la suite $r=0$, 1 ou 2:
les chemins 0-colori\'es sont les chemins de Dyck,
les chemins 1-colori\'es sont les chemins de Motzkin ordinaires,
et les chemins 2-colori\'es seront dits {\em bicolori\'es}.

Un {\em pr\'efixe} de Motzkin est un chemin
de Motzkin \'eventuellement incomplet,
c'est \`a dire qui ne se termine pas
forc\'ement sur l'axe des abscisses;
la {\em hauteur} d'un pr\'efixe d\'esigne l'ordonn\'ee
du point d'arriv\'ee.

La {\em longueur} d'un chemin compte le nombre de pas;
un chemin peut \^etre vide: il est alors de longueur nulle.
En notant un pas ascendant $a$, un pas descendant $b$,
et un pas horizontal de couleur $i$ , $c_i$  $(1 \leq i \leq r)$\/,
on associe \`a un chemin de  Motzkin $r$-colori\'e
un {\em mot de Motzkin}, sur un alphabet \`a $r+2$ lettres;
les mots de Dyck correspondent au cas $r=0$;
le pas $c_1$ sera parfois not\'e simplement $c$.

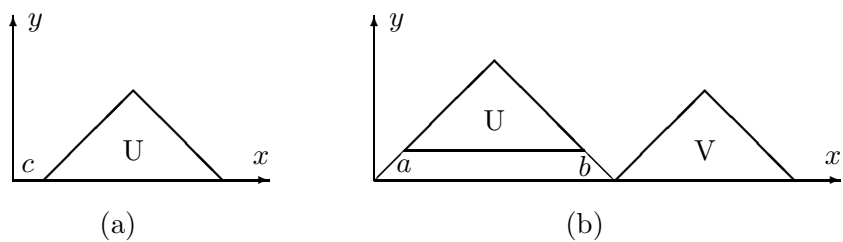
\begin{figure}[htbp]
\begin{centering}
\setlength{\unitlength}{0.02cm}
\begin{picture}(550,140)(0,-30)


\thinlines
\put(0,0){\vector(1,0){170}}
\put(165,15){\makebox(0,0){$x$}}
\put(0,0){\vector(0,1){110}}
\put(15,105){\makebox(0,0){$y$}}
\put(70,-30){\makebox(0,0){(a)}}

\put(0,0){\line(1,0){20}}
\put(10,10){\makebox(0,0){$c$}}

\thicklines
\put(20,0){\line(1,1){60}}
\put(80,60){\line(1,-1){60}}
\put(80,20){\makebox(0,0){U}}


\thinlines
\put(240,0){\vector(1,0){310}}
\put(545,15){\makebox(0,0){$x$}}
\put(240,0){\vector(0,1){110}}
\put(255,105){\makebox(0,0){$y$}}
\put(380,-30){\makebox(0,0){(b)}}

\put(240,0){\line(1,1){20}}
\put(380,20){\line(1,-1){20}}
\put(260,10){\makebox(0,0){$a$}}
\put(380,10){\makebox(0,0){$b$}}

\thicklines
\put(260,20){\line(1,1){60}}
\put(320,80){\line(1,-1){60}}
\put(400,0){\line(1,1){60}}
\put(460,60){\line(1,-1){60}}

\put(260,20){\line(1,0){120}}

\put(320,40){\makebox(0,0){U}}
\put(460,20){\makebox(0,0){V}}

\end{picture}
\caption{D\'ecompositions de chemins}
\label{fig.chemins}
\end{centering}
\end{figure}

Un chemin de Motzkin non vide $W$ se d\'ecompose selon le sch\'ema de la
\fig{chemins}, en $W=c_i U$ -- cas (a) --
ou en $W=a U b V$ -- cas (b) ;
$U$ et $V$ d\'esignent \`a nouveau
des chemins de Motzkin (\'eventuellement vides):
on se trouve dans le cas (a) ou (b) selon que W d\'ebute par
un pas horizontal ou ascendant; dans ce dernier cas,
$V$ d\'ebute d\`es le {\em premier}
retour sur l'axe des abscisses, d'o\`u
l'unicit\'e de la \decomp.

En d\'esignant par $\Mz$ le langage des mots de Motzkin $r$-colori\'es,
on en d\'eduit l'\'equation:
\begin{equation}
  \label{eq.motzkin}
  \Mz = 1 + \sum_{1 \leq i \leq r} c_i \Mz + a \Mz b \Mz
\end{equation}
(o\`u la somme $\sum c_i \Mz$ dispara\^\i t pour $r=0$\/).
La \senum\ des mots de Motzkin selon leur longueur v\'erifie donc:
\[
  \mu(t) = 1 + r t \mu(t) + t^2 \mu(t)^2
\]
dont les solutions sont \eq{sol-eqr-strict} et \eq{sol-eqr}
pour $r=1$ et $r=2$, avec 
$\lambda(t) = t \mu(t)$.

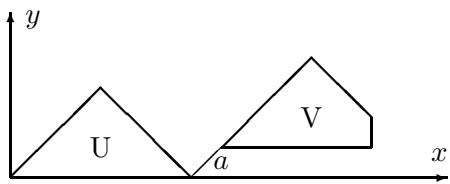
\begin{figure}[htbp]
\begin{centering}
\setlength{\unitlength}{0.02cm}
\begin{picture}(290,110)(240,0)

\thinlines
\put(240,0){\vector(1,0){290}}
\put(525,15){\makebox(0,0){$x$}}
\put(240,0){\vector(0,1){110}}
\put(255,105){\makebox(0,0){$y$}}

\put(360,0){\line(1,1){20}}
\put(380,10){\makebox(0,0){$a$}}

\thicklines
\put(240,0){\line(1,1){60}}
\put(300,60){\line(1,-1){60}}
\put(380,20){\line(1,1){60}}
\put(440,80){\line(1,-1){40}}

\put(380,20){\line(1,0){100}}
\put(480,20){\line(0,1){20}}

\put(300,20){\makebox(0,0){U}}
\put(440,40){\makebox(0,0){V}}

\end{picture}
\caption{D\'ecomposition d'un pr\'efixe}
\label{fig.schema-prefixe}
\end{centering}
\end{figure}

Un pr\'efixe $W$ d'un chemin de Motzkin,
ou bien est un chemin de Motzkin complet,
ou bien se d\'ecompose selon la \fig{schema-prefixe},
en un chemin de Motzkin $U$ suivi d'un pr\'efixe $V$,
soit $W=UaV$.
Cette \decomp\ est unique,
car $U$ se termine au {\em dernier} passage
sur l'axe des abscisses. Donc la \senum\ $\rho(t)$
des pr\'efixes est li\'ee \`a celle des chemins complets par:
\[
  \rho(t) = \mu(t) [1+t\rho(t)]
\]

\label{celibat1}
Nous appellerons {\em c\'elibataire} un pas ascendant
conduisant du niveau $y$ au niveau $y+1$
--~ie. un pas $(x,y) \rightarrow (x+1,y+1)$~--,
et qui n'est suivi d'aucun pas descendant ramenant le chemin
au niveau $y$. Si on it\`ere la d\'ecomposition ci-dessus
d'un pr\'efixe $W$, on obtient:
\begin{equation}
  \label{eq.prefixe}
  W = U_0 \: a \: U_1 \: a \ldots U_l
\end{equation}
o\`u les occurrences de $a$ d\'esignent les pas ascendants c\'elibataires,
et les $U_i$ des mots de Motzkin (\'eventuellement vides).
Cette \decomp\ est unique, et la hauteur du chemin est \'egale
au nombre $l$ de pas c\'elibataires.

L'analogie entre les \decomp s de chemins et d'animaux,
fournit une bijection
r\'esum\'ee dans le \theor\ suivant,
d\^u \`a \cite{gouyou} (mais avec une construction totalement
diff\'erente):

\begin{theoreme}
  \label{th.animaux-chemins}
  Il existe une bijection entre, d'une part:
  animaux dirig\'es, de source ponctuelle, sur r\'eseau carr\'e
  (resp. triangulaire),
  de taille $n$, et d'autre part:
  pr\'efixes de Motzkin (resp. pr\'efixes bicolori\'es),
  de longueur $n-1$.
  Dans cette bijection, les animaux de \dld\ $k$
  correspondent aux chemins de hauteur~$k$;
  en particulier les \eqr s correspondent aux chemins de Motzkin.
\end{theoreme}

\preuve la bijection se construit de mani\`ere r\'ecursive;
d\'etaillons-la dans le sens chemin $\rightarrow$ animal,
et notons-la $\beta$;
rappelons que $|U|$ d\'esigne la longueur d'un chemin $U$,
et $|P|$ la taille d'un animal $P$.
Au chemin vide (de longueur nulle) on associe
l'animal r\'eduit \`a une cellule (de taille 1).

Soit $W$ un chemin de Motzkin non vide; si $W=c_i U$,
on associe r\'ecur\-sivement \`a $U$ un animal en \'equerre
$M=\beta(U)$, et on pose $\beta(W)=L=xM$,
avec les notations de \eq{eqr1};
si $W$ est un chemin de Motzkin ordinaire, la base de
$M$ est translat\'ee horizontalement de -1 par rapport \`a $x$,
base de $L$,
comme sur la \fig{equerres};
si $W$ est un chemin de Motzkin bicolori\'e, la base de
$M$ est translat\'ee de -1 ou de 0 selon que $W$
d\'ebute par $c_1$ ou $c_2$: on obtient alors un animal
sur r\'eseau triangulaire.

De m\^eme, si $W=aUbV$,
on associe r\'ecursivement \`a $U$ et $V$
des \eqr s $M=\beta(U)$ et $N=\beta(V)$,
et on pose $\beta(W)=L=xMN$,
avec les notations de \eq{eqr2}. La base de
$M$ est translat\'ee horizontalement de -1
par rapport \`a $x$, base de $L$,
et $N$ a sa base sur la m\^eme fibre que $x$.

Si $W$ est un pr\'efixe de Motzkin qui ne retourne pas sur
l'axe des abscisses, alors $W=UaV$, o\`u $a$ d\'esigne un pas
ascendant c\'elibataire;
on associe r\'ecursivement \`a $U$ et $V$
une \'equerre $L=\beta(U)$ et un animal $Q=\beta(V)$,
et on pose $\beta(W)=P=LQ$,
avec les notations de \eq{pyreq}. La base de
$Q$ est translat\'ee horizontalement de +1 par rappport \`a celle de $L$.

On v\'erifie que ces constructions conservent la relation
$|\beta(W)|=|W|+1$
et que la \dld\ de $\beta(W)$ s'accro\^\i t en m\^eme temps que
la hauteur de $W$ (rappelons que la hauteur de $W$ est d\'efinie comme
l'ordonn\'ee du point d'arriv\'ee)\cqfd

\medskip

Cette bijection est \'ecrite en langage C dans la \sct{code},
et la \sct{illustrations} contient des exemples
(avec $c_1$ not\'e $c$, et $c_2$ not\'e $d$).
Le \thref{animaux-chemins} confirme les valeurs trouv\'ees
pour les coefficients des s\'eries \eq{sol-eqr}
et \eq{sol-pyr}, gr\^ace \`a la proposition suivante:

\begin{proposition}
Le nombre de mots de Motzkin bicolori\'es de longueur $n-1$
est \'egal au nombre de mots de Dyck de longueur $2n$\/,
donn\'e par le nombre de Catalan:
\[
  C_n = \frac{1}{n+1}\binome{2n}{n}
      = \frac{1}{2n+1}\binome{2n+1}{n}\ .
\]
Le nombre de pr\'efixes de Motzkin bicolori\'es de longueur $n-1$
est \'egal au nombre de pr\'efixes de Dyck de longueur $2n-1$\/,
donn\'e par le coefficient binomial:
\[
  \binome{2n-1}{n} = \frac{1}{2}\binome{2n}{n}\ .
\]
\end{proposition}

\preuve cette proposition est classique en combinatoire;
la premi\`ere partie d\'ecoule simplement de l'identit\'e entre \sgens:
\[
  [1 + \lambda(t) ]^2 = 1 + 2 \lambda(t) + \lambda(t)^2 \ .
\]

Nous en donnons une preuve bijective pour \^etre complets.
Pour la pre\-mi\`ere affirmation, il suffit d'associer \`a
un mot de Motzkin bicolori\'e $W$,
un mot de Dyck
{\em non vide} $W'$ d\'efini de fa\c con r\'ecursive
-- $\epsilon$ d\'esigne le chemin vide:
\[
  \epsilon \rightarrow ab\ , \hspace{1em}
  c_1 U \rightarrow abU'\ , \hspace{1em}
  c_2 U \rightarrow aU'b\ , \hspace{1em}
  aUbV \rightarrow aU'bV'\ .
\]
C'est une bijection, car l'\'equation $D=1+aDbD$, v\'erifi\'ee
par la s\'erie \carac\ des mots de Dyck --~voir \eq{motzkin}~--,
peut s'\'ecrire: $$D-1=ab+ab(D-1)+a(D-1)b+a(D-1)b(D-1)\ .$$

Pour la seconde partie de la proposition,
on part de $  W = U_0 \: a \: U_1 \: a \ldots U_l $
--~voir \eq{prefixe}~--,
on associe \`a chaque mot de Motzkin bicolori\'e $U_i$ un mot de Dyck $U_i'$
comme ci-dessus, on appelle $U_i''$ le mot $U_i'$
priv\'e de son pas descendant final $b$, et on associe \`a W
le mot: $$W'' = U_0'' \: a \: U_1'' \: a \ldots U_l''\ . $$
C'est une bijection car on peut retrouver la \decomp\
ci-dessus de $W''$ en rep\'erant ses pas ascendants c\'elibataires
de rang pair. On v\'erifie aussi que $|W''| = 2|W|+1$.

Enfin, un pr\'efixe de Dyck de longueur impaire
$2n-1$ poss\`ede un nombre impair $2k-1$ de
pas ascendants c\'elibataires; 
en rempla\c cant les $k$ premiers
par des pas descendants, on obtient un chemin
joignant l'origine au point $(2n-1, -1)$. C'est encore une bijection,
car les pas descendants ainsi introduits sont les pas
\label{celibat2}
{\em descendants c\'elibataires}\/: un pas descendant
conduisant du niveau $y+1$ au niveau $y$
--~ie. un pas $(x,y+1) \rightarrow (x+1,y)$~--,
est appel\'e c\'elibataire
s'il n'est {\em pr\'ec\'ed\'e} d'aucun pas ascendant
conduisant du niveau $y$ au niveau $y+1$,
autrement dit s'il est le premier atteignant le niveau $y$.
Ceci explique que les pr\'efixes de Dyck soient compt\'es
par un coefficient binomial\cqfd

\medskip

\begin{figure}[htbp]
\begin{centering}
\setlength{\unitlength}{0.02cm}
\newcommand{\chapeau}[1]{%
  \thicklines
  \begin{picture}(80,40)
  \put(0,0){\line(1,1){40}}
  \put(0,0){\line(1,0){80}}
  \put(40,40){\line(1,-1){40}}
  \put(40,15){\makebox(0,0){$U_{#1}$}}
  \end{picture}
}

\newcommand{\ascendant}{%
  \thinlines
  \begin{picture}(20,20)
  \put(0,0){\line(1,1){20}}
  \put(20,5){\makebox(0,0){$a$}}
  \end{picture}
}

\newcommand{\descendant}{%
  \thinlines
  \begin{picture}(20,0)
  \put(0,0){\line(1,-1){20}}
  \put(0,-15){\makebox(0,0){$b$}}
  \end{picture}
}

\begin{picture}(620,160)(0,-80)

\thinlines
\put(0,0){\vector(1,0){620}}
\put(615,15){\makebox(0,0){$x$}}
\put(0,-80){\vector(0,1){160}}
\put(15,75){\makebox(0,0){$y$}}

\put(0,0){\chapeau{0}}
\put(80,0){\descendant}
\put(100,-20){\chapeau{1}}
\put(180,-20){\descendant}
\multiput(210,-50)(10,-10){3}{\circle*{2}}

\put(240,-80){\chapeau{k}}
\put(320,-80){\ascendant}
\put(340,-60){\chapeau{k+1}}
\put(420,-60){\ascendant}
\multiput(450,-30)(10,10){5}{\circle*{2}}
\put(500,20){\chapeau{k+l}}

\end{picture}
\caption{Factorisation de Catalan d'un chemin quelconque}
\label{fig.catalan}
\end{centering}
\end{figure}
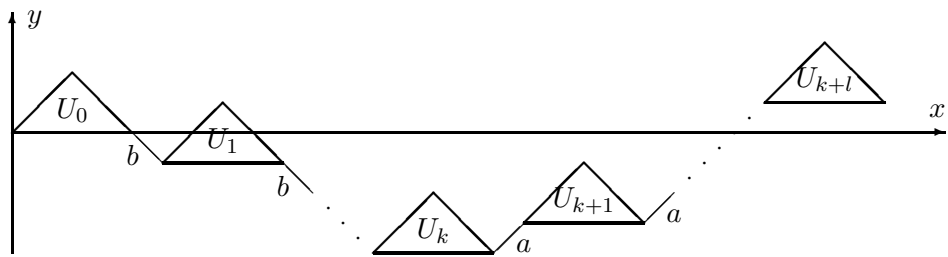

La technique qui consiste \`a rep\'erer, dans un chemin
(ou mot) quelconque $W$, les pas c\'elibataires,
conduit \`a la {\em factorisation de Catalan}\/:
\begin{equation}
  \label{eq.catalan}
  W = U_0 \: b \: U_1 \: b \ldots U_k
  \: a \: U_{k+1} \: a  \ldots U_{k+l}
\end{equation}
qui comporte $k$ (resp. $l$) pas c\'elibataires descendants
(resp. ascendants), et est unique: voir \fig{catalan}.

Les pr\'efixes de Motzkin sont les chemins sans pas
c\'elibataire descendant ($k=0$), et le nombre $l$ de pas
c\'elibataires ascendants est \'egal \`a la hauteur du chemin;
les chemins de Motzkin
sont ceux qui ne poss\`edent aucun pas c\'elibataire,
ascendant ou descendant ($k=l=0$).

A un pr\'efixe de Motzkin de hauteur $l$, on peut associer $l+1$
mots en rempla\c cant les $j\ (0 \leq j \leq l)$ premiers pas
c\'elibataires ascendants par des pas descendants; on en d\'eduit,
comme dans \cite{gouyou}, que la moyenne des hauteurs des
pr\'efixes de Motzkin $r$-colori\'es de longueur $n$ vaut:
\[
  \frac{(r+2)^n}{m_n} - 1
\]
o\`u $m_n$ d\'esigne le nombre de ces pr\'efixes;
comme la largeur moyenne d'un animal est \'evidemment
le double de la \dld\ moyenne,
le  \thref{animaux-chemins} fournit le corollaire suivant,
emprunt\'e \`a \cite{gouyou}:

\begin{corollaire}
  Soit $a_n$ le nombre d'animaux de source ponctuelle,
  sur r\'eseau carr\'e ou triangulaire,
  de taille $n$ (donn\'e par le \thref{enum-animaux});
  la largeur moyenne de ces animaux vaut:
  \[
  2 \:\frac{(r+2)^{n-1}}{a_n} - 2
  \]
avec $r=1$ sur r\'eseau carr\'e, $r=2$ sur r\'eseau triangulaire\cqfd
\end{corollaire}

Appelons {\em compacte} une configuration, sur un r\'eseau, aussi
dense que possible (en respectant bien entendu la r\`egle d'exclusion
entre sites voisins); sur un r\'eseau lin\'eaire, la distance
entre deux cellules d'une configuration compacte est donc \'egale
\`a deux. Un animal de source compacte correspond \`a un empilement
de base compacte.
La factorisation de Catalan fournit imm\'ediatement une bijection entre mots
et animaux de source compacte: \`a chaque facteur
$U_i$ de \eq{catalan}, on associe par le \thref{animaux-chemins} un animal
en \eqr\ $L_i$, et \`a $W$ on associe le produit des \eqr s,
en espa\c cant de {\em deux} unit\'es les bases de $L_0 \ldots L_k$,
puis seulement d'une unit\'e les bases de $L_k \ldots L_l$.
Voir une illustration \sct{illustrations}.
Nous obtenons ainsi une nouvelle preuve de
l'\'etonnant th\'eor\`eme d\'ecouvert par  \cite{gouyou}:

\begin{corollaire}
  \label{th.compact}
  Le nombre d'animaux de taille $n$ et de source compacte
  sur r\'eseau carr\'e (resp. triangulaire) est \'egal \`a $3^{n-1}$
  (resp. $4^{n-1}$)\cqfd
\end{corollaire}

\section{Animaux plans \alea s}
\label{sec.aleatoires}

\subsection{Code}
\label{sec.code}

Les th\'eor\`emes \ref{th.animaux-chemins} et \ref{th.compact}
de la \sct{motzkin} fournissent des bijections explicites
entre animaux et mots sur un alphabet \`a $r+2$ lettres
($r=1$ pour un r\'eseau carr\'e,
 $r=2$ pour un r\'eseau triangulaire); ils permettent
de construire facilement
un programme de g\'en\'eration \alea\ d'animaux plans,
avec distribution uniforme.
La complexit\'e de l'algorithme est lin\'eaire en moyenne:
le temps d'ex\'ecution est proportionnel, en moyenne,
\`a la taille $n$ de l'animal.
La m\'emoire utilis\'ee est proportionnelle \`a $n$,
et on ne peut \'evidemment faire mieux si on m\'emorise l'animal;
sinon on pourrait se contenter d'une taille m\'emoire
proportionnelle, en moyenne, \`a $\sqrt{n}$\/.
Nous donnons explicitement ce programme, sous forme de
fonctions \'ecrites en langage C.

\begin{figure}[htbp]
\begin{quote}
\small
\begin{verbatim}
static int nb_tirages; char Mot  [TAILLE_MAX];

static char hasard (n)
int n; {return (char) ( 'a' + lrand48() % n );}

static void mot (n, r)
int n, r; {int i;
  for (i = 0; i < n; i++) Mot [i] = hasard (r + 2);
  Mot [n] = '\0';
}

static void prefixe (n, r)
int n, r; {int i, h; char pas;
  for (nb_tirages = i = h = 0; i < n; nb_tirages++) {
    Mot [i++] = pas = hasard (r + 2);
    if (pas == 'a') h++;
    else if (pas == 'b')
      if (--h < 0) i = h = 0;  /* on recommence tout */ 
  }
  Mot [n] = '\0';
}
\end{verbatim}
\end{quote}
\caption{Code pour la g\'en\'eration \alea\ de pr\'efixes de Motzkin}
\label{fig.code-motzkin}
\end{figure}

Il faut commencer
par construire un mot au hasard,
de longueur $n$, sur un alphabet de $r+2$ lettres.
Si l'on garde n'importe quel mot
(fourni par une suite de tirages de lettres, effectu\'es
\`a l'aide d'un g\'en\'erateur de nombres \alea s
de la biblioth\`eque num\'erique Unix standard),
l'animal correspondant a une source compacte (voir \thref{compact}).
Si l'on souhaite construire un animal de source ponctuelle,
il ne faut garder que les pr\'efixes de Motzkin (voir \thref{animaux-chemins});
dans ce cas on recommence la proc\'edure depuis le d\'ebut,
d\`es que le chemin associ\'e traverse l'axe des abscisses;
cet algorithme tr\`es simple ne d\'etruit pas le caract\`ere
uniforme de la distribution, et il est montr\'e dans
\cite{barcucci} que sa complexit\'e
est en moyenne $2n$ -- voir aussi \cite{alonso}.

La fonction {\em prefixe (n, r)}, qui calcule dans la variable
globale {\em Mot}, un pr\'efixe de Motzkin \alea\ de longueur $n$,
sur un alphabet \`a $r+2$ lettres, est donn\'ee \fig{code-motzkin}.
La variable globale {\em nb\_tirages} permet de calculer
la complexit\'e (exp\'erimentale) de l'algorithme, et de v\'erifier,
si on le souhaite,
les pr\'evisions de \cite{barcucci}. Nous donnons aussi la
fonction {\em mot~(n,~r)}, qui fournit un mot \alea\ quelconque,
pour \^etre complets.
Le g\'en\'erateur de nombres \alea s utilis\'e,
{\em lrand48~()}, doit en principe \^etre initialis\'e
avant emploi, par appel \`a {\em srand48~()}: voir le manuel Unix.

\begin{figure}[htbp]
\begin{quote}
\small
\begin{verbatim}
static celibataires_ascendants (n)
int n; {int i, h = 0, hmin = 0; char pas;
  for (i = n-1; i >= 0; i--) {
    pas = Mot [i];
    if (pas == 'b') h++;
    else if (pas == 'a')
      if (--h < hmin) { hmin = h; Mot [i] = 'A'; }
  }
}

static celibataires_descendants (n)
int n; {int i, h = 0, hmin = 0; char pas;
  for (i = 0; i < n; i++) {
    pas = Mot [i];
    if (pas == 'a') h++;
    else if (pas == 'b')
      if (--h < hmin) { hmin = h; Mot [i] = 'B'; }
  }
}
\end{verbatim}
\end{quote}
\caption{Code pour marquer les c\'elibataires}
\label{fig.code-celibataires}
\end{figure}

La construction de la bijection entre mots et animaux,
fait jouer un r\^ole tout \`a fait particulier aux pas
c\'elibataires (voir d\'efinitions pages \pageref{celibat1}
et \pageref{celibat2}). Ici 
on marque les pas c\'elibataires en transformant la lettre
correspondante en majuscule.
Nous incluons ces deux fonctions \'el\'ementaires,
pour \^etre explicites, \fig{code-celibataires}.

\begin{figure}[htbp]
\begin{quote}
\small
\begin{verbatim}
static int
  compteur,            /* compte les lettres  */
  fibre [FIBRE_MAX];   /* hauteurs des fibres */

static void equerre (i)
int i; {
  empiler (i);
  switch (Mot [compteur++]) {
    case 'a':
      equerre (i - 1); equerre (i); break;
    case 'c':
      equerre (i - 1); break;
    case 'd':
      equerre (i); break;
  }
}

void animal (n, triangulaire, compact)
int n, triangulaire, compact; {int i, r;
  
  if (triangulaire) r = 2; else r = 1;
  for (i = 0; i < FIBRE_MAX; i++) fibre [i] = -1;

  if (compact) {
    mot (n - 1, r);
    celibataires_descendants (n - 1);
  }
  else
    prefixe (n - 1, r);
  celibataires_ascendants (n - 1);

  for (i = compteur = 0; compteur < n; i++) {
    equerre (i); if (Mot [compteur - 1] == 'B') i++;
  }
}
\end{verbatim}
\end{quote}
\caption{Code pour transformer un mot en animal}
\label{fig.code-animal}
\end{figure}

La fonction {\em equerre~(i)}, sur la \fig{code-animal},
construit un animal en \eqr\ dont la base a pour abscisse~$i$,
autrement dit est situ\'ee sur la fibre~$i$.
Les hauteurs des fibres sont conserv\'ees dans la variable {\em fibre},
et la fonction {\em empiler~(i)}, d\'etaill\'ee plus loin,
ajoute une cellule sur la fibre num\'ero $i$ de l'animal.
Le reste du code de la fonction {\em equerre} est une traduction
directe de la preuve du \thref{animaux-chemins}; les pas $c$ et $d$
sont les deux types de pas horizontaux, not\'es $c_1$ et $c_2$
\sct{motzkin}; le pas $d$ n'est utilis\'e que sur r\'eseau
triangulaire. On notera l'absence de r\'ef\'erence explicite
au pas $b$, qui sert de s\'eparateur entre deux \eqr s:
apr\`es lecture d'un pas ascendant $a$, et construction r\'ecursive
d'une premi\`ere \eqr, par l'appel {\em equerre~(i-1)},
on trouve forc\'ement un pas descendant $b$,
sinon le pas ascendant serait c\'elibataire, et not\'e $A$;
on appelle alors {\em equerre~(i)}.

La fonction {\em animal (n, triangulaire, compact)} est la fonction
principale, qui construit un animal de taille~$n$,
sur r\'eseau carr\'e ou triangulaire,
selon la valeur du param\`etre bool\'een {\em triangulaire};
la source est ponctuelle ou compacte
selon la valeur du param\`etre bool\'een {\em compact}.
Dans cette fonction, on initialise la taille de l'alphabet
selon la nature du r\'eseau, on rase les fibres,
et on construit un mot \alea\ quelconque,
ou un pr\'efixe de Motzkin, selon que la source de l'animal
est compacte ou ponctuelle. Apr\`es transformation des
pas c\'elibataires (un pas ascendant est not\'e $a$ ou $A$
selon qu'il est ordinaire ou c\'elibataire; idem pour un
pas descendant, not\'e $b$ ou $B$), qui s\'eparent les \eqr s,
la fonction construit
de fa\c con r\'ep\'etitive des \eqr s, d\'ecal\'ees l'une par
rapport \`a l'autre de +1 ou de +2, selon que la lettre qui
s\'epare deux \eqr s d\'enote un pas c\'elibataire
ascendant ou descendant:
voir la fin de la preuve du \thref{animaux-chemins},
et la preuve du \thref{compact}.
 La variable globale {\em compteur},
incr\'ement\'ee \`a chaque lecture d'une lettre du mot,
permet d'en rep\'erer la fin.

\begin{figure}[htbp]
\begin{quote}
\small
\begin{verbatim}
typedef struct {int x, y;} Cellule;
Cellule Animal [TAILLE_MAX];

static void empiler (i)
int i; {int h = -2, j = i + DECALAGE, jmax, k;
  for (k = j - 1; k <= j+1; k++)
    if (h < fibre [k]) {h = fibre [k]; jmax = k;}
  h++; if (jmax == j) h++;
  Animal [compteur].x = i;
  Animal [compteur].y = fibre [j] = h;
}
\end{verbatim}
\end{quote}
\caption{Code pour empiler une cellule}
\label{fig.code-empiler}
\end{figure}

La fonction {\em empiler~(i)}, sur la \fig{code-empiler},
est \'el\'ementaire:
elle ajoute une cellule sur la fibre~$i$ de l'animal,
en fonction des tailles des fibres voisines.
Les cellules sont stock\'ees dans le tableau
{\em Animal}; la constante {\small \tt DECALAGE}
sert \`a convertir les abscisses des fibres, qui peuvent
\^etre n\'egatives, en indices de tableaux.
Si la fibre voisine la plus haute, de hauteur~$h$,
est la fibre~$i$ elle-m\^eme,
le r\'eseau est triangulaire, et l'ordonn\'ee de la nouvelle
cellule vaut $h+2$, au lieu de $h+1$ en g\'en\'eral.

\subsection{Illustrations}
\label{sec.illustrations}

La fonction {\em animal}, de la \sct{code}, nous a servi
a construire les illustrations de cette section.
L'un des auteurs a construit un logiciel d'affichage
pour {\em X-Window}, en utilisant les outils {\em Xview},
pour {\em OpenLook}.

\begin{figure}
  \includegraphics[height=10cm]{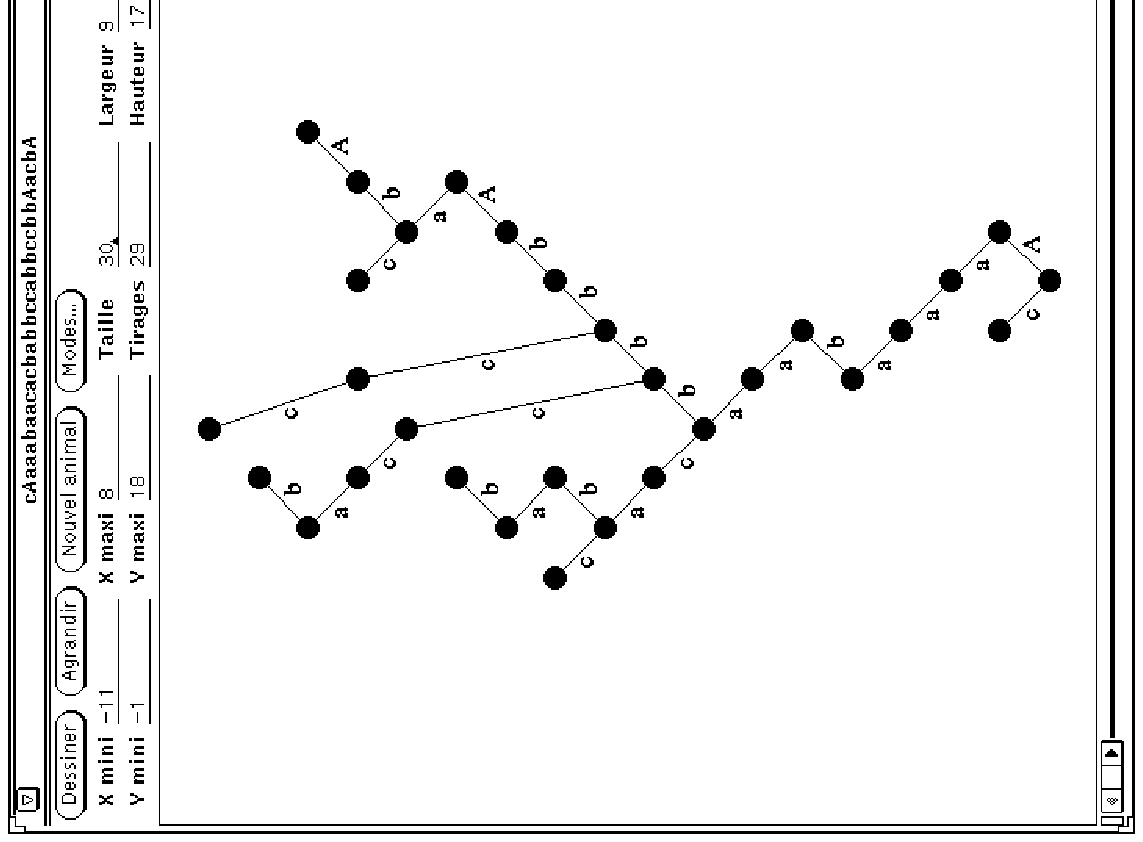}
  \includegraphics[height=10cm]{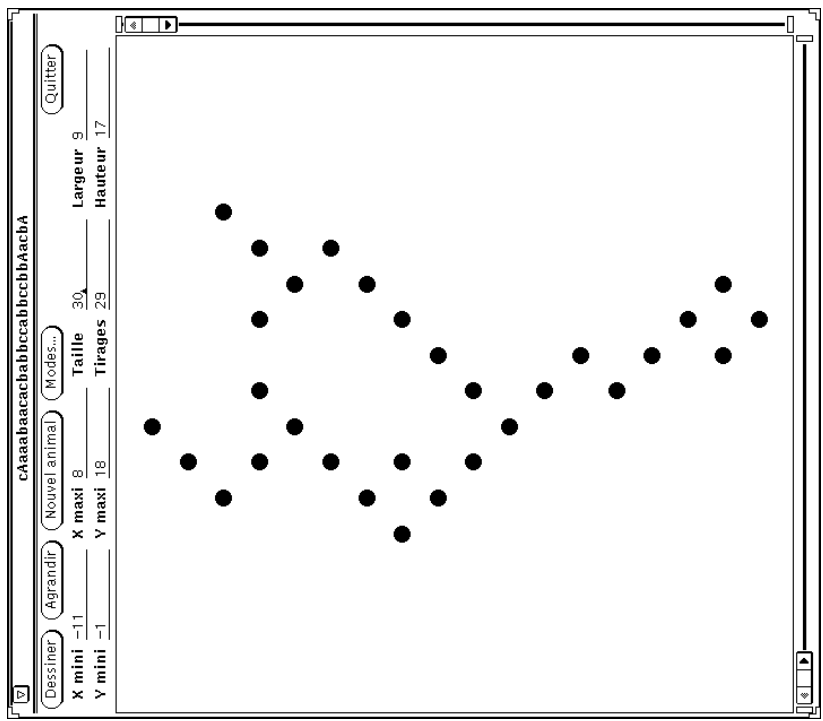}
  \caption{Exemple de construction d'un animal sur r\'eseau carr\'e}
  \label{fig.carre30}
\end{figure}

\begin{figure}
  \includegraphics[height=10cm]{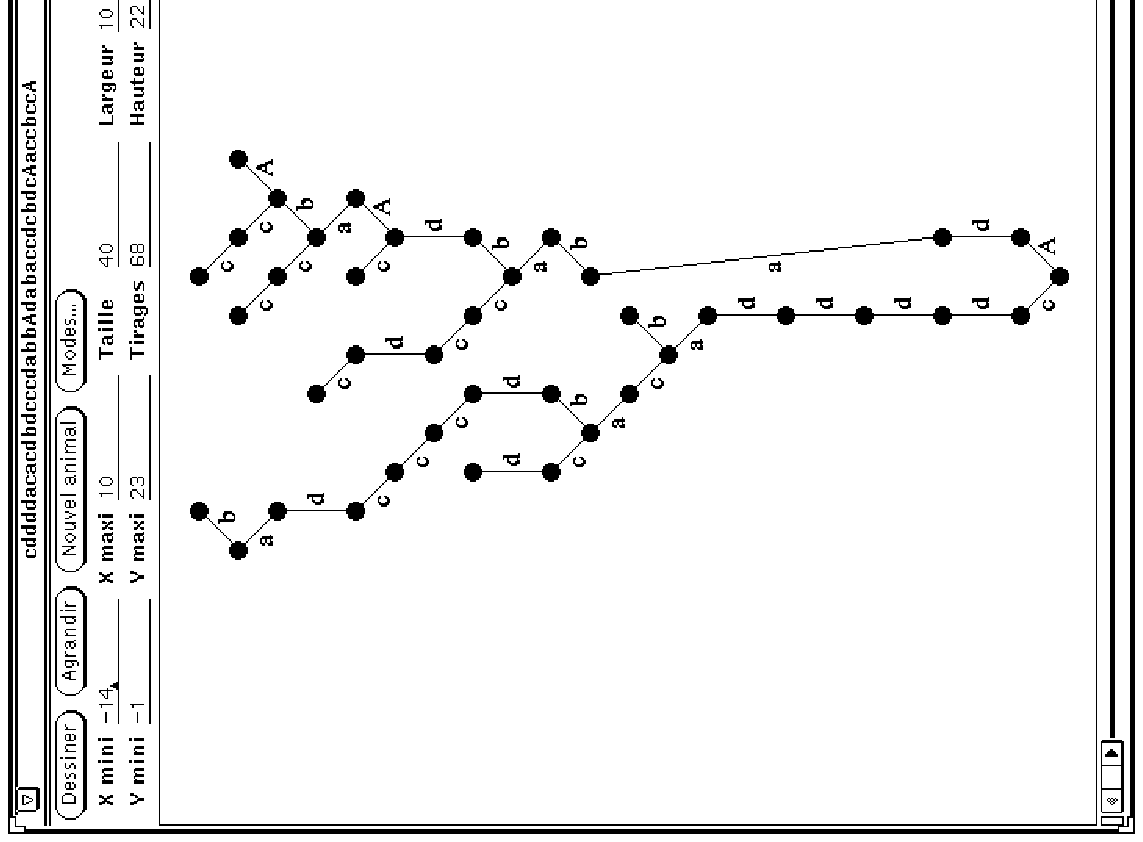}
  \includegraphics[height=10cm]{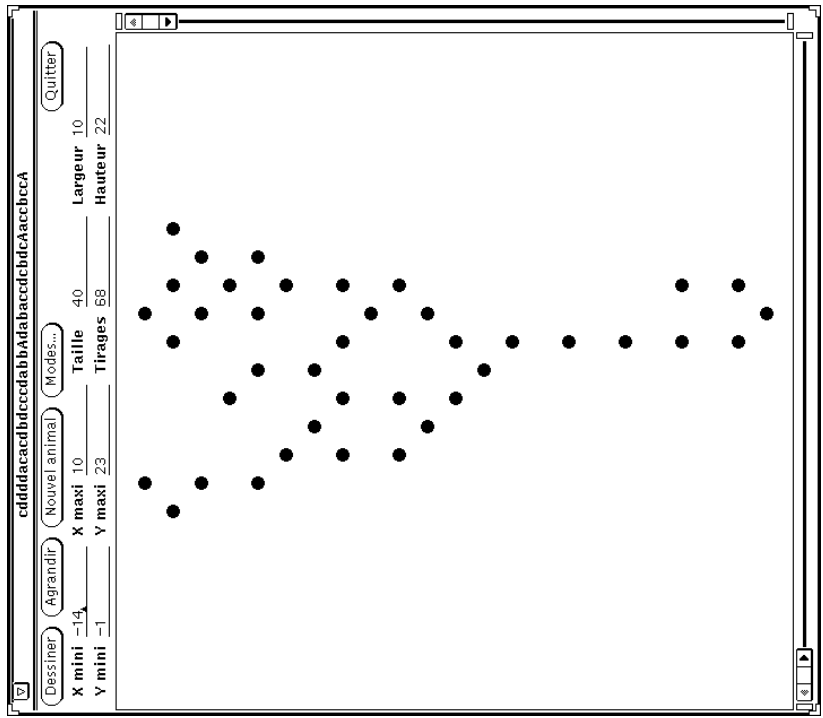}
  \caption{Exemple de construction d'un animal sur r\'eseau triangulaire}
  \label{fig.triangulaire}
\end{figure}

\begin{figure}
  \includegraphics[height=10cm]{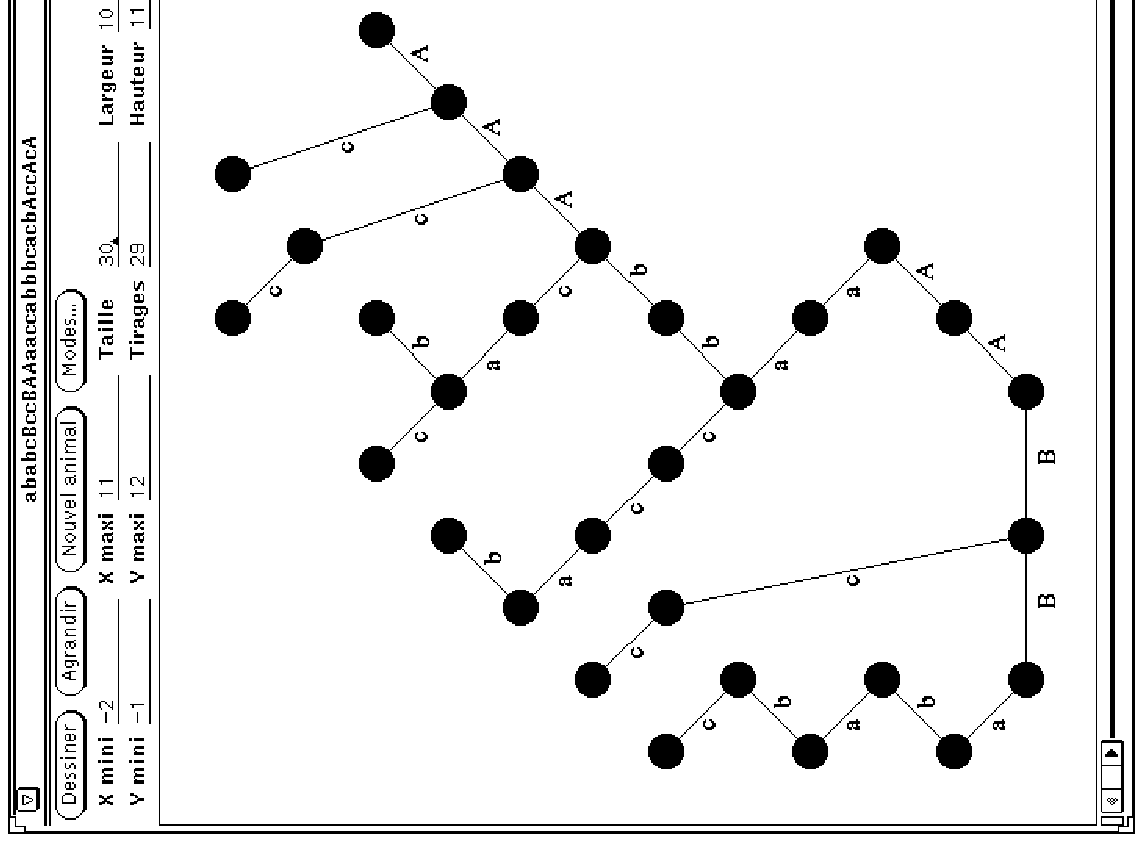}
  \includegraphics[height=10cm]{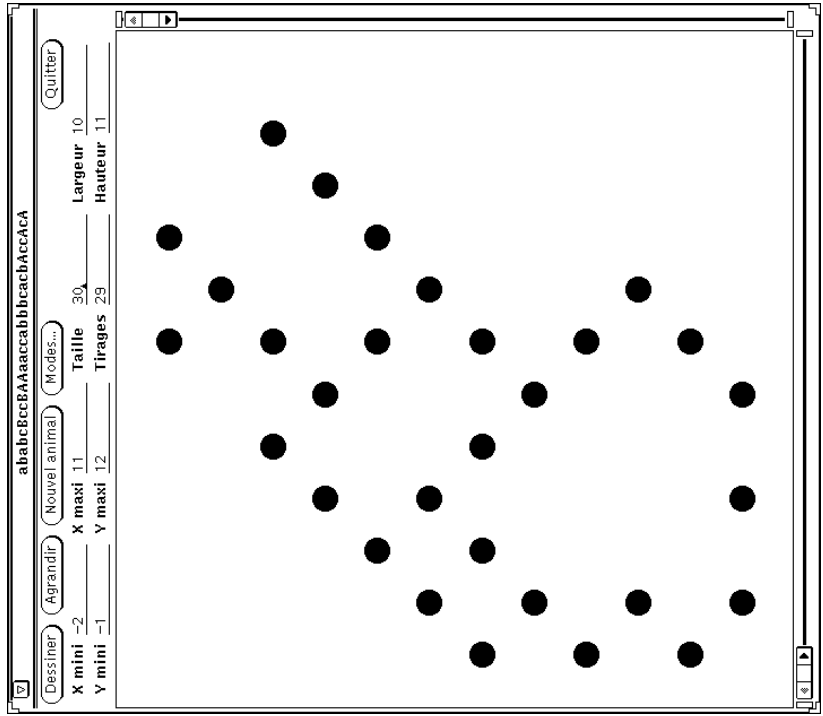}
  \caption{Exemple de construction d'un animal de source compacte}
  \label{fig.compact30}
\end{figure}

Chacune des figures \ref{fig.carre30},
\ref{fig.triangulaire}, et \ref{fig.compact30},
contient sur sa partie droite un arbre,
appel\'e {\em arbre de guingois}
par les auteurs dans~\cite{penaud}.
Le parcours en profondeur de cet arbre,
de gauche \`a droite,
donne le mot associ\'e \`a l'animal,
qui se trouve en titre en haut de la fen\^etre,
avec les pas c\'elibataires not\'es $A$ ou $B$;
cet arbre r\'esume la construction de l'animal \`a partir du mot.

Comme l'algorithme de construction d'un animal
ajoute des \eqr s de gauche \`a droite,
on peut les lire de droite \`a gauche,
\`a condition, chaque fois qu'on en lit une (sous-arbre
issu du point le plus bas de la fibre la plus \`a droite),
de l'enlever avant de lire la suivante.

\begin{figure}
  \includegraphics[width=\textwidth]{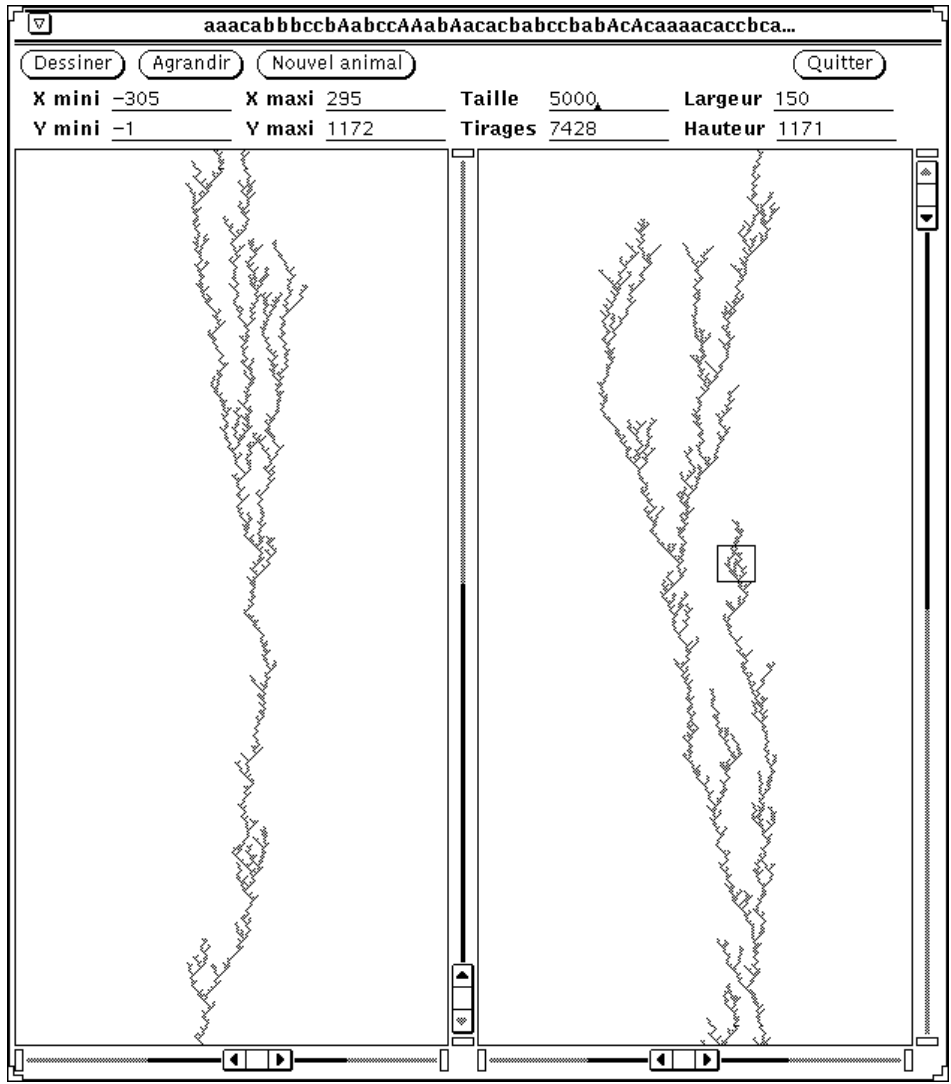}
  \caption{Exemple d'animal de grande taille sur r\'eseau carr\'e}
  \label{fig.carre5000}
\end{figure}

\begin{figure}
  \includegraphics[width=\textwidth]{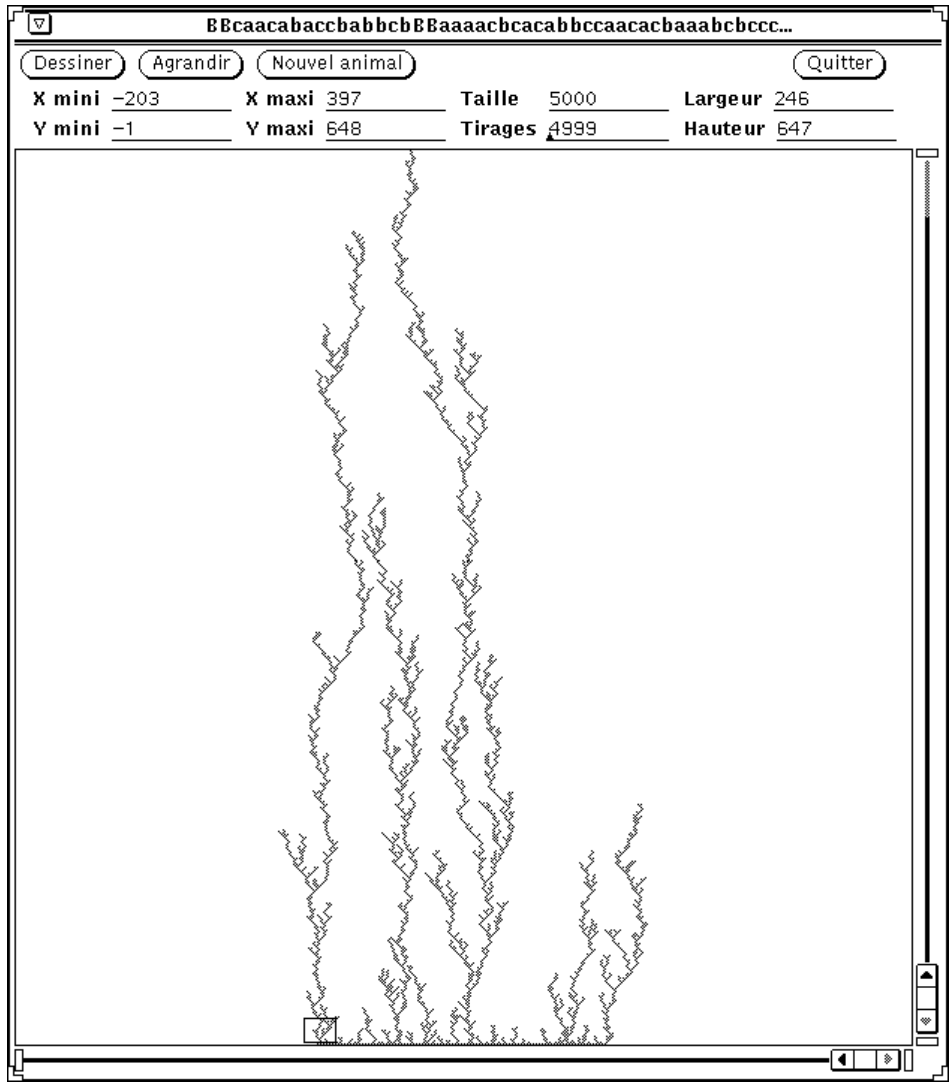}
  \caption{Exemple d'animal de grande taille et de source compacte}
  \label{fig.compact5000}
\end{figure}

\begin{figure}
  \includegraphics[height=10cm]{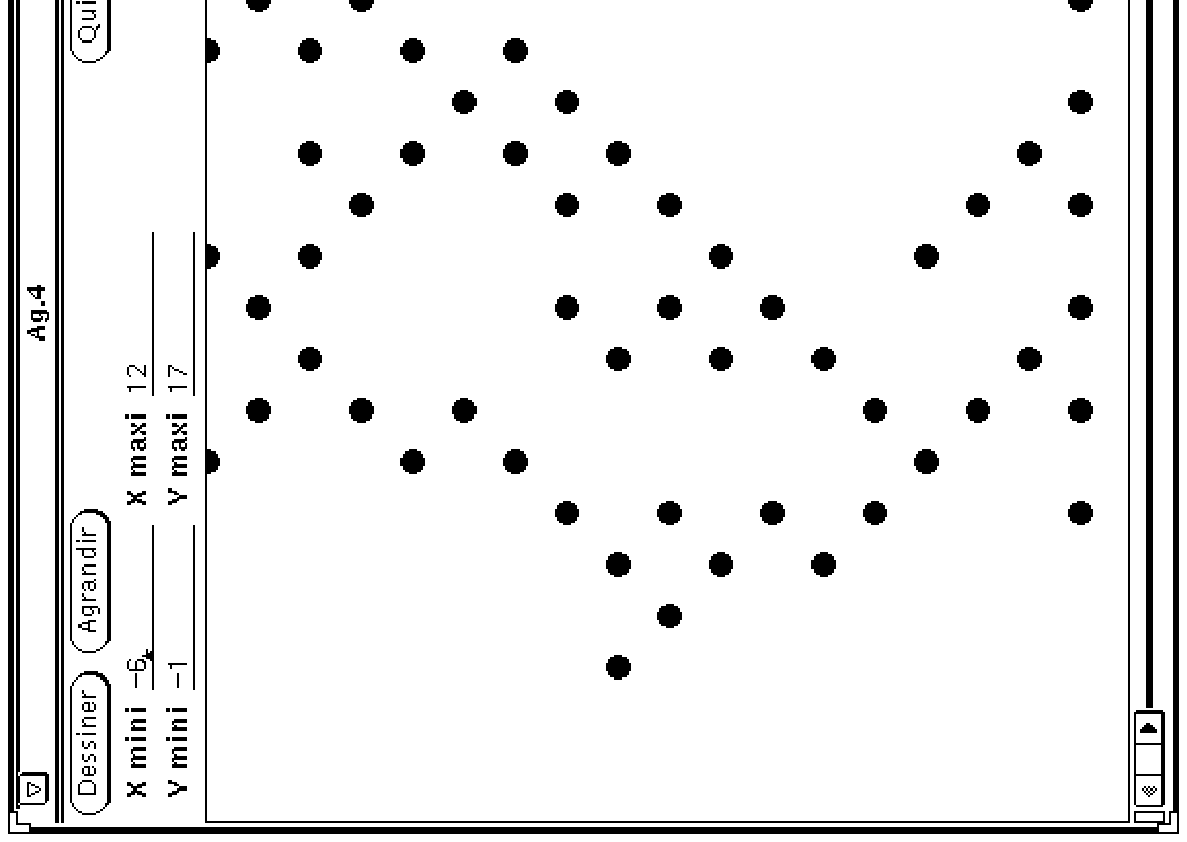}
  \includegraphics[height=10cm]{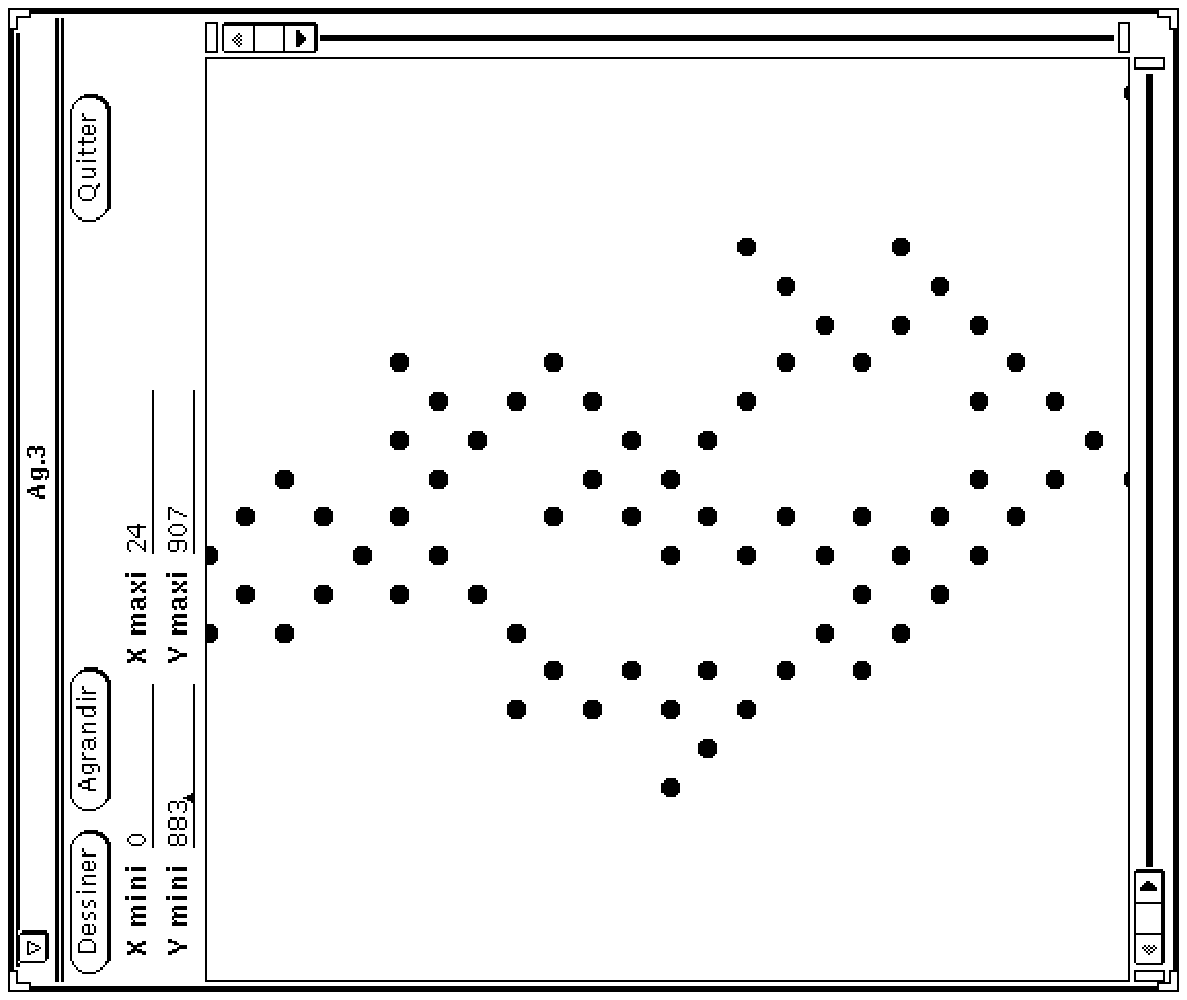}
  \caption{Agrandissements de r\'egions des figures
    \protect\ref{fig.carre5000} et \protect\ref{fig.compact5000}}
  \label{fig.agrandi}
\end{figure}

Les figures \ref{fig.carre5000} et \ref{fig.compact5000}
pr\'esentent des animaux de grande taille (5000 cellules),
qui sont obtenus sur une station Sun-4 en un temps de
l'ordre de la seconde, affichage compris.
La \fig{carre5000} utilise la possibilit\'e de scinder
une fen\^etre en deux:
la moiti\'e gauche pr\'esente
la partie inf\'erieure de l'animal,
et la moiti\'e droite la partie sup\'erieure,
comme on peut le voir \`a la position des ``ascenseurs'' verticaux.
Le logiciel permet
d'agrandir \`a volont\'e une partie de l'animal
(d\'elimit\'ee par un rectangle): la \fig{agrandi}
illustre ces possibilit\'es. Le logiciel affiche aussi, pour
chaque animal, sa taille,
sa largeur (mesur\'ee sur l'axe des abscisses),
sa hauteur (mesur\'ee sur l'axe des ordonn\'ees),
et le nombre de tirages
n\'ecessaire pour obtenir le mot \alea\ associ\'e
(ce dernier renseignement n'est utile que pour les animaux de
source ponctuelle). L'utilisateur peut modifier la taille
avant de demander la g\'en\'eration d'un nouvel animal.

Les animaux montr\'es ici ont \'evidemment \'et\'es produits
de fa\c con \alea, mais les auteurs ont s\'electionn\'e
des exemples qui leur ont sembl\'e instructifs ou esth\'etiques;
l'aspect moyen d'un animal est souvent plus simple.






\section{Perspectives}
\label{sec.perspectives}

Cet article a pr\'esent\'e un calcul sur les animaux dirig\'es,
qui permet de d\'ecouvrir une structure cach\'ee de ces derniers,
symbolis\'ee par les arbres des figures \ref{fig.carre30},
\ref{fig.triangulaire}, et \ref{fig.compact30}.
Parmi les probl\`emes ouverts concernant les animaux plans,
le principal est l'\'evaluation de leur hauteur moyenne
(appel\'ee aussi longueur).
Plus g\'en\'eralement, on devrait pouvoir \'etudier plus finement
la structure de ces animaux, et obtenir par exemple
leur dimension fractale, ou le nombre moyen
de composantes connexes du compl\'ementaire.

Le probl\`eme majeur reste celui des mod\`eles de gaz sur
des r\'eseaux en dimension $d>1$, qui correspondent \`a certains
animaux dirig\'es en dimension $d+1$; le cas du r\'eseau hexagonal
a \'et\'e r\'esolu par R.Baxter
(solution du mod\`ele des {\em heagones durs}, voir \cite{baxter}),
et la s\'erie g\'en\'eratrice des pyramides sur un tel r\'eseau
est alg\'ebrique, mais une interpr\'etation combinatoire des
\'equations obtenues reste \`a trouver. Le cas du r\'eseau carr\'e
(mod\`ele des {\em carr\'es durs})
reste enti\`erement ouvert.

\newcommand{\auteur}[1]{{\bf #1}}
\newcommand{\collection}[1]{\ \linebreak[2] \rm #1.}
\newcommand{\editeur}[1]{\ \linebreak[2] {\rm #1.}}
\newcommand{\titre}[1]{{\em #1}}

\end{document}